\documentclass[12pt,a4paper]{article}%
\usepackage{amsmath}
\usepackage{enumitem}
\usepackage{graphicx}
\usepackage{epsfig}%
\usepackage{amsfonts}%
\usepackage{amssymb}
\usepackage{makeidx}         
\usepackage{color}
\usepackage{bm}
\usepackage{soul}
\usepackage{appendix}
\usepackage{subfig}
\RequirePackage{mathrsfs}
\usepackage{yfonts}
\usepackage{appendix}
\usepackage{url}

\usepackage{psfrag}
\usepackage{fancyhdr}
\usepackage{epsfig}
\usepackage{algorithm}
\usepackage{algorithmic}
\usepackage{amsmath,amssymb,graphicx}
\usepackage{bm}
\usepackage{cancel}
\usepackage{amsfonts}
\usepackage{pgfplots} 
\usepackage{graphicx}
\input{pgfart.sty}
\usepackage{algorithmic}

\RequirePackage{mathrsfs}

\DeclareMathAlphabet{\Bi}{OT1}{cmm}{b}{it}
\providecommand{\F}{\mathfrak}

\newcommand{\refeq}[1]{Eq.~(\ref{#1})}

\newtheorem{prop}{Proposition}[section]

\newtheorem{algo}{Algorithm}[section]
\newtheorem{defi}{Definition}[section]

\newtheorem{rem}{Remark}[section]

\makeatletter\@addtoreset{equation}{section}\makeatother
\makeatletter\@addtoreset{figure}{section}\makeatother
\makeatletter\@addtoreset{table}{section}\makeatother
\newtheorem{theorem}{Theorem}[section]

\newtheorem{ex}[theorem]{Example}
\newtheorem{lemma}[theorem]{Lemma}

\newenvironment{proof}[1][Proof]{\noindent \emph{#1.} }{\hfill \ 
\rule{0.5em}{0.5em}}

\newcommand{\vepsilon}{\varepsilon}

\newcommand{\vTheta}{\varTheta}

\definecolor{myred}{rgb}{1, 0.2, 0.2}

\newtheorem{notation}{Notation}[section]

\def\c{\textnormal{\textbf{c}}}
\def\u{\textnormal{\textbf{u}}}
\def\c{\textnormal{\textbf{c}}}
\def\x{\textnormal{\textbf{x}}}

\def\f{\textnormal{\textbf{f}}}

\def\w{\omega}
\def\wo{\mathop \omega \limits^{\circ}}
\def\w1{\omega_1}
\def\w2{\omega_2}

\def\Th{\mathcal{T}_h}

\def\TH{\mathcal{T}_H}
\def\G{\mathcal{G}}
\def\tA{\tilde{A}}
\def\F{\mathcal{F}}
\def\O{\mathcal{O}}
\def\H{\mathcal{H}}    
\def\Ht{\tilde{H}}    
\def\F{\mathcal{F}}    
\def\RR{\mathbb{R}}

\def\l{\lambda}

\def\beginproof{\textbf{Proof: }}
\def\endproof{\hfill$\blacksquare$\\}

\def\Proof{\textbf{Proof: }}
\DeclareMathAlphabet{\mathitbf}{OML}{cmm}{b}{it}
\def\jmath{j}

\newcommand{\di}{\mathrm{d}}

\def\RR{\mathbb{R}}

\def\H{\mathcal{H}}

\def\u{\bm{u}}
\def\bZ{{Z}}
\def\pp{\pi}

\def\f{\mathbf{f}}

\begin{document}

\title{Partial inversion of the elliptic operator to speed up computation of likelihood in Bayesian inference}

\author{Alexander Litvinenko}



%
\maketitle
\tableofcontents
\newcommand{\slugmaster}{%
\slugger{juq}{xxxx}{xx}{x}{x--x}}

\pagestyle{myheadings}
\thispagestyle{plain}
\markboth{A. Litvinenko}{Partial inversion of elliptic operator to speed up computation of Bayesian inference}

\renewcommand{\thefootnote}{\fnsymbol{footnote}}

\footnotetext[2]{E-mail: litvinenko@uq.rwth-aachen.de. RWTH Aachen, Aachen, Germany}
\renewcommand{\thefootnote}{\arabic{footnote}}

\newpage
\begin{abstract}
Often, when solving forward, inverse or data assimilation problems, only a part of the solution is needed. As a model, we consider the stationary diffusion problem.
We demonstrate an algorithm that can compute only a part or a functional of the solution, without calculating the full inversion operator and the complete solution.
It is a well-known fact about partial differential equations that the solution at each discretisation point depends on the solutions at all other discretisation points.
Therefore, it is impossible to compute the solution only at one point, without calculating the solution at all other points. 
The standard numerical methods like a conjugate gradient or Gauss elimination compute the whole solution and/or the complete inverse operator.
We suggest a method which can compute the solution of the given partial differential equation 1) at a point; 2) at few points; 3) on an interface; or a functional of the solution, without computing the solution at all points. The required storage cost and computational resources will be lower as in the standard approach.

With this new method, we can speed up, for instance, computation of the innovation in filtering or the likelihood distribution, which measures the data misfit (mismatch). Further, we can speed up
the solution of the regression, Bayesian inversion, data assimilation, and Kalman filter update problems.

Applying additionally the hierarchical matrix approximation, we reduce the cubic computational cost to almost linear $\mathcal{O}(k^2n \log^2 n)$, where $k\ll n$ and $n$ is the number of degrees of freedom. 

Up to the hierarchical matrix approximation error, the computed solution is exact.
One of the disadvantages of this method is the need to modify the existing deterministic solver.
\end{abstract}


\textbf{AMS}
65F10,  
60H15,  
60H35,  
65C30  

\begin{table}[htbp!]
\caption{Notation}
\begin{tabular}{|c|l|}
\hline 
HDD & suggested here the hierarchical domain decomposition method \\ \hline
$u\vert_{\gamma}$ & restriction of the solution $u$ onto the interface $\gamma$\\ \hline
$h$, $H$            &grid step sizes on fine and coarse meshes \\ \hline
$\Omega$, $\partial\Omega$ & computational domain and its boundary \\ \hline
$Z$ & random parameter vector $Z=(Z_1,...,Z_{n_Z})$\\ \hline
$\vTheta$  & space where parameter $\bZ=(Z_1,...,Z_{n_{z}})$ is defined \\ \hline
$\omega$, $\partial\omega$ & local subdomain and its boundary \\ \hline
$V_h$, $V_H$    & two finite element spaces, $V_H\subset V_h$ \\ \hline
$\f$, $\f_h$, $\f_H$            & the right hand side, discretized on fine ($h$) and coarse ($H$) meshes \\ \hline
$\u$, $\u_h$, $\u_H$         & the solution, computed on fine  ($h$) and coarse ($H$) meshes \\ \hline
$\kappa(x,Z)=e^{q(x,Z)}$   & uncertain permeability coefficient, depends on parameter vector $Z$\\ \hline
$\mathcal{V}_N$ & vector space spanned on the basis $\{\varphi_1(x),\ldots,\varphi_N(x)\}$ \\ \hline
$I$, $I_N$ & index sets \\ \hline \hline
$\Th$, $\TH$ & fine and coarse triangulations \\ \hline
$T_{\Th}$  & hierarchical domain decomposition tree\\ \hline
$\gamma_{\omega}$ & interface in the domain $\omega \subset \Omega$ (also call ``internal" boundary)\\ \hline
$\Gamma_{\omega}=\partial\omega$ & boundary (also call ``external''  boundary)\\ \hline
$d_{\omega}$ & $d_{\omega}: = \left( {\left( {f_i} \right)_{i \in I(\omega)} ,\left( {g_i }
\right)_{i \in I(\partial_{\omega} )} } \right)=(f_{\omega},g_{\omega})$ \\
& a composed vector consisting of the right-hand side restricted to $\omega$  \\ 
& and the Dirichlet boundary values $g_{\omega}=u_h\vert_{\partial\omega}$ \\ \hline
$\F_h$, $\G_h$ & two operators, such that $u_h=\F_hf_h+\G_hg_h$,   \\ \hline  
$\hat{y}$ & true observations \\ \hline
$y= \hat{y}+ \vepsilon$ & noisy observations \\ \hline
$\Phi^g_{\omega}:\RR^{I(\partial\omega)}\rightarrow \RR^{I(\gamma_{\omega})}$ & maps the boundary data defined on $\partial\omega$ to the \\
&data defined on the interface $\gamma_{\omega}$  \\ \hline
 $\Phi^f_{\omega}:\RR^{I(\omega)}\rightarrow \RR^{I(\gamma_{\omega})}$ & maps the right-hand side data defined on $\omega$ to the \\
 &data defined on $\gamma_{\omega}$.\\ \hline
$\Psi^f_{\omega}:\RR^{ I(\omega) }\rightarrow  \RR^{ I(\partial\omega) }$ &  maps the whole subdomain to the external boundary \\ \hline
 $\Psi^g_{\omega}:\RR^{ I(\partial\omega) }\rightarrow  \RR^{ I(\partial\omega) }$ & maps the external boundary to the external boundary \\ \hline
pdf & probability density function \\ \hline
PCG & preconditioned conjugate gradient \\ \hline
\end{tabular}
\label{tab:notations}
\end{table}%
\section{Introduction}
\label{sec:intro}

We further develop the method, initially introduced in \cite{MYPHD,Hackbusch2012, FlorianDiss, HDDdoc} and in Chapter 12 of \cite{HackHMBook}.
With this method, we will be able to compute the solution in a subdomain, in a point, the mean value over a subdomain and other functionals $F(u)$ without computing the full inverse operator and the complete solution.  
Similar ideas were considered in \cite{babb2016accelerated, martinsson2015hierarchical}.
This method can be very practical for speeding up the solution of the inverse and data assimilation problems, which appear in many science and engineering applications such as weather prediction, oil recovery, and subsurface flow.
Under the inverse problem, we understand the estimation of unknown model parameters from (noisy) measurements. Under the data assimilation problem, we understand improving the existing mathematical model when the new measurement data become available. 

The forward problem we consider is the diffusion problem with uncertain or unknown diffusion coefficient. A typical task is not only to solve the forward problem but also to identify this unknown coefficient. 
Initially, some prior probability density function for the diffusion coefficient is assumed.  
Then the Bayesian inference is applied to update this density.
The Bayesian inference is a statistical inference method in which Bayes' theorem is used to update the probability for a hypothesis as more evidence or information becomes available.

Computing the likelihood function in the Bayesian formula requires multiple solutions of the forward problem and could be time-consuming.
Depending on the available measurements, the complete solution of the (forward) diffusion problem could be unnecessary, {rather} only a part of the solution or a functional of the solution is needed. 
Computing only a part of the solution will make the whole computing process faster and less time-consuming.

Typically, the available measurement data is a functional $F(u)$ of the solution
$u$. The misfit (or mismatch) function is the difference between the simulated data and the measurement values \cite{rosic2012sampling,rosic2013parameter,hermann2016inverse,matthies2016parameter,pajonk2012deterministic,rosic2011direct}.
Below we will show how to simulate these measurement data directly without computing the whole solution. Particularly, we will show that calculating the full inverse operator is unnecessary.

{One possible application of our method is the data driven research, a very popular topic nowadays. In this research the available datasets are used either to improve (enrich) the existing mathematical model (often a system of PDEs), or to discover the governing system of PDEs. Another example when fast calculation of a part of the solution is required, is computing the mean square error when comparing the training and computed datasets.
In \cite{Raissi2019}, authors design data-driven algorithms for inferring solutions to various partial differential equations.
They introduce neural networks that are trained to solve supervised learning tasks while respecting any given laws of physics described by general nonlinear PDEs. Their goal is to solve two classes of problems: data-driven solution and data-driven discovery of PDEs.}

The structure of this paper is the following. In Section~\ref{sec:intro} we give our motivation by introducing
the stochastic forward problem and the Bayesian updating procedure for computing posterior density function of the uncertain diffusion coefficient.
The main ingredient and the main contribution --- the hierarchical domain decomposition (HDD) method --- is contained in Section~\ref{sec:HDD}. Details of the HDD method, including two algorithms ``Leaves to Root'' and ``Root to Leaves'', are shown in Section~\ref{sec:HDDConstrProcess}. The hierarchical (denoted by $\mathcal{H}$) - matrix technique to speed up the HDD method
is explained in Section~\ref{sec:Hmatrices}. 
Section \ref{sec:functional} 
explains how to use the HDD method to compute a functional of the solution without computing the whole solution. 
Particularly, it explains how to compute the mean value in a small subdomain. The novelty here is that the whole solution is not available, only a small part of it. 
In the last section, we conclude the main achievements.

\begin{figure}[htbp]
\centerline{ 
\includegraphics[width=0.24\textwidth]{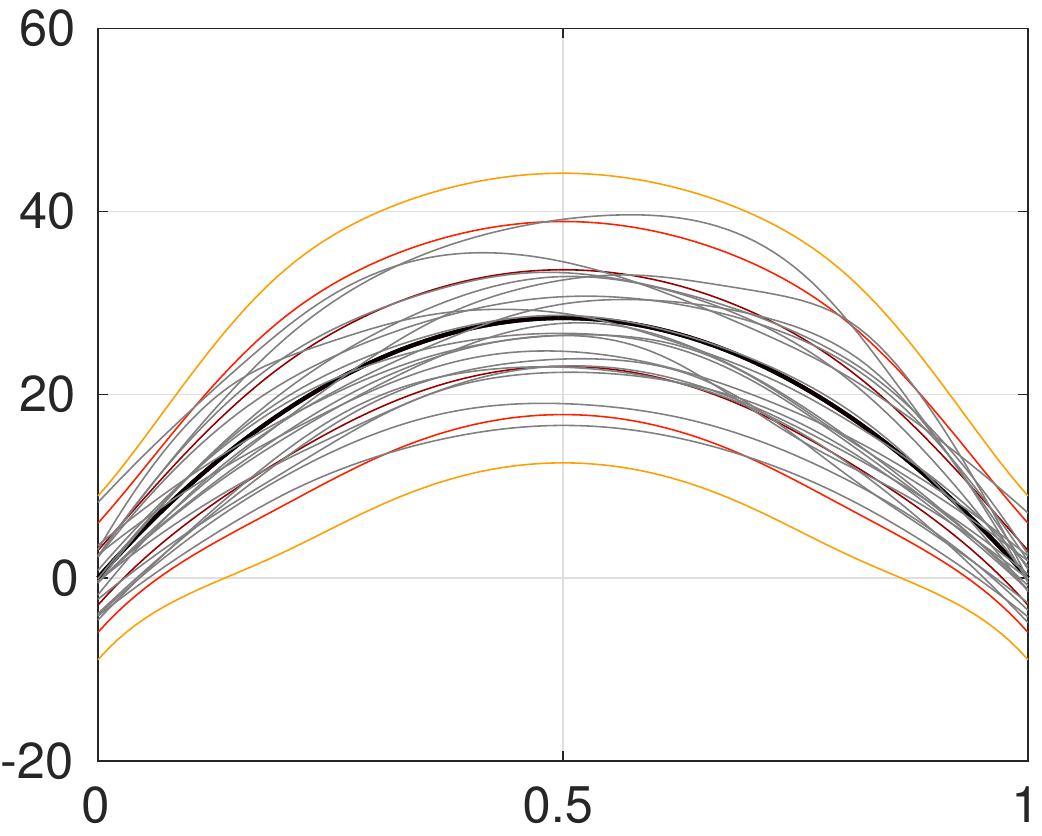}
 \includegraphics[width=0.24\textwidth]{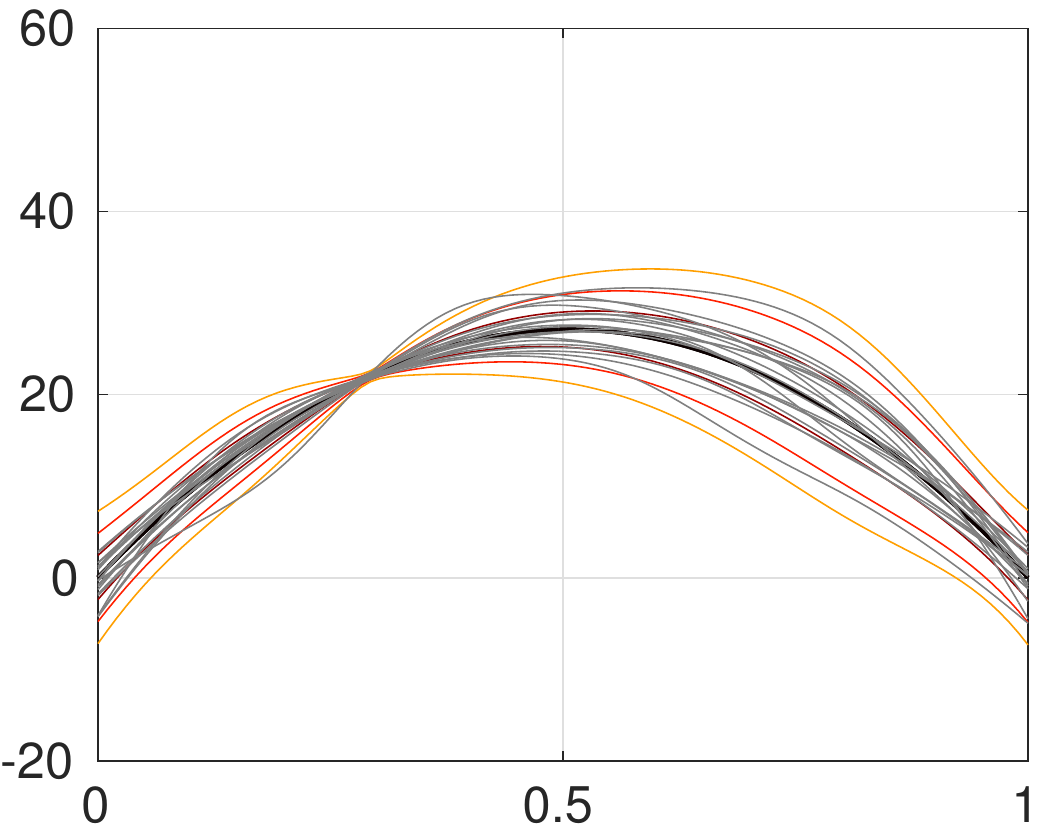}
 \includegraphics[width=0.24\textwidth]{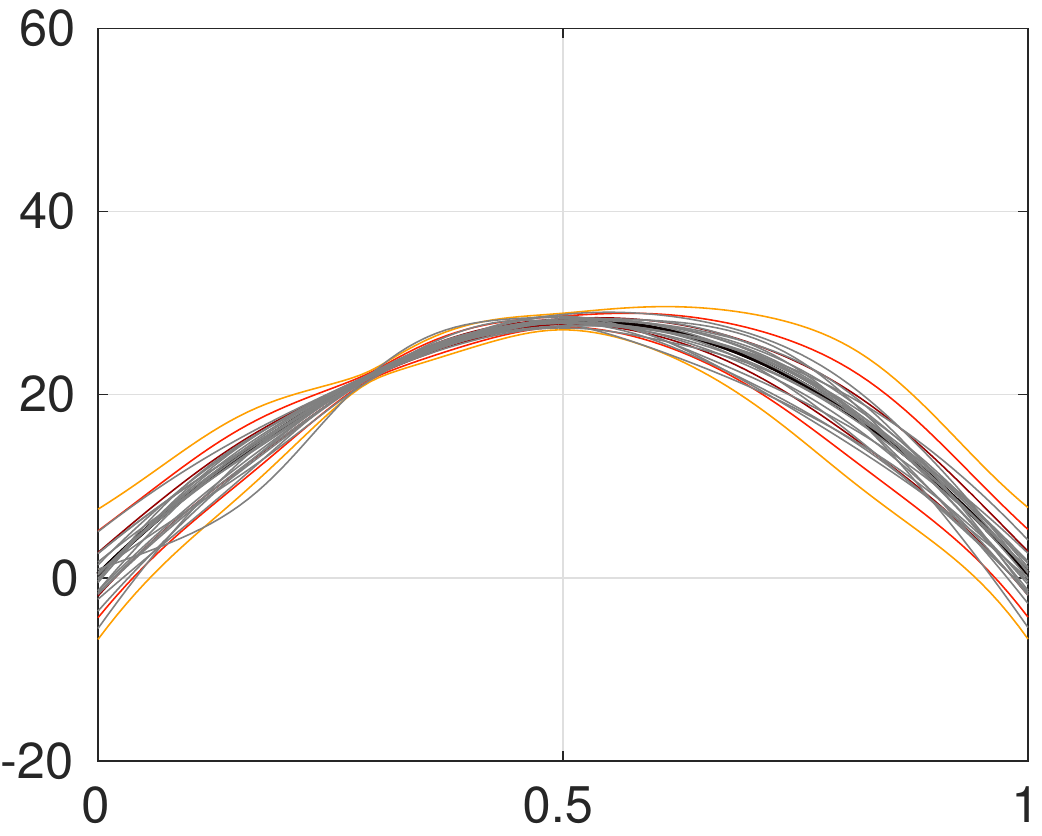}
 \includegraphics[width=0.24\textwidth]{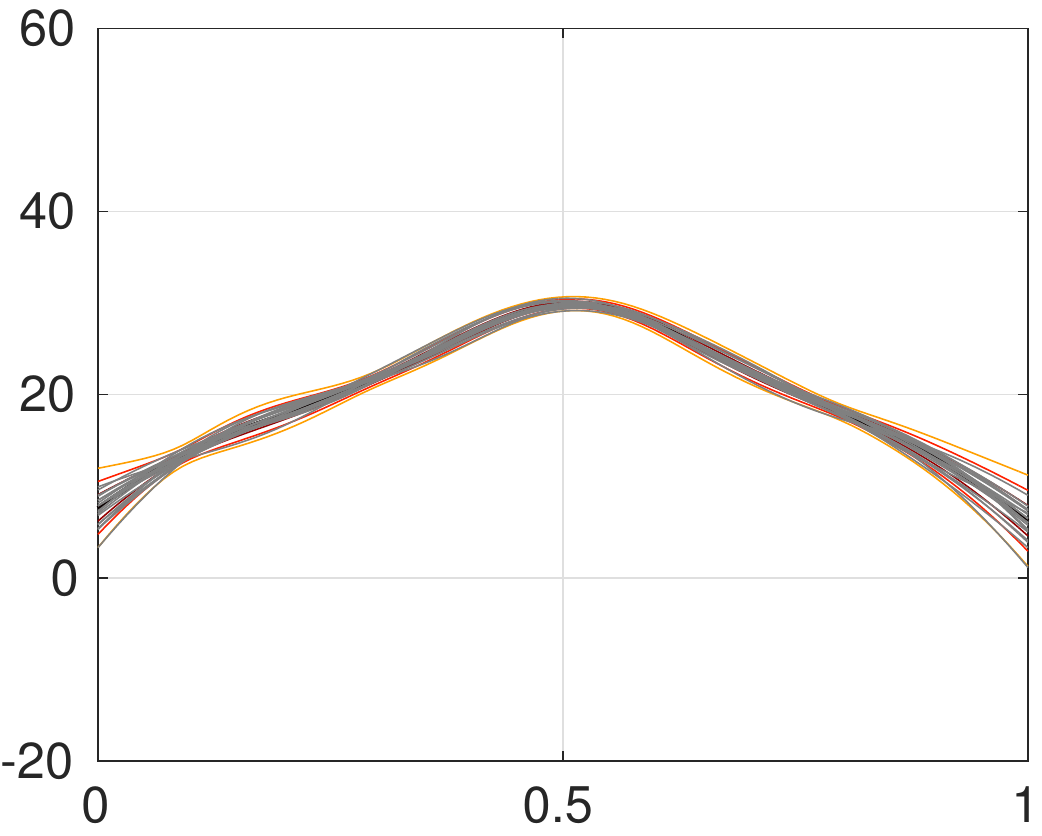}}
\caption{(left) 20 original realisations of the solution $u$ (2nd, 3rd, 4th) the same realisations after update; the mean value (bold line) and $± {\displaystyle \pm }{1,2,3}$ standard deviations (red, orange and yellow lines). }
\label{fig:update_u}
\end{figure}

\textbf{Example.}
This example shows how the solution (or measurements) in only a few points can reduce the uncertainty.
Consider an elliptic PDE with uncertain coefficient and the right hand side as in Eq.~\ref{eq:elliptic}, but in 1D, on the interval $[0,1]$.
We pose uncertain Dirichlet boundary conditions  $g(0,\xi)$ and $g(1,\xi)$, where $\xi$ is a Gaussian random variable.
Assume three measurements at locations $x=\{0.3,0.5,0.8\}$ are given. The mean values $\overline{u}(0.3)=22$, $\overline{u}(0.5)=28$, $\overline{u}(0.8)=18$ and the standard deviations are $\{0.2,0.3,0.3\}$ respectively.
The following computations are done with the stochastic Galerkin library \textit{sglib}, written by E. Zander at TU Braunschweig \footnote{\url{https://github.com/ezander/sglib}}.
In Fig.~\ref{fig:update_u} twenty realisations of the uncertain solution $u(x)$ before and after an update are shown. The mean value  (dark bold line) and $± {\displaystyle \pm }\{1,2,3\}$ standard deviations (red, orange and yellow lines) are shown. The left picture shows realisations, obtained with some prior assumption about distribution of random diffusion coefficient $\kappa$. In the following three pictures of the updated solutions are shown, after taking into account one, two and three measurements.
To conclude this example, it is practical to have a numerical method, which can efficiently compute a part of the solution or a solution in a few points, without computing the complete solution.
%
\subsection{Main idea}
\label{sec:MainIdea}

The main ingredients of the developed approach are the weak formulation, the hierarchical (or recursive) domain decomposition technique, the finite element method, and the Schur complement. Additionally, to speed up matrix operations and reduce the overall storage cost, we approximate the involved operators and the 
Schur complement in the hierarchical matrix format \cite{Part1,GH03,HackHMBook}.

The \textbf{novelty} of this work is the application of the HDD method for faster computation of the innovation in filtering or the likelihood distribution, which measures the data misfit (mismatch).

The forward problem we consider is an elliptic boundary value problem with uncertain $L^{\infty}$ coefficients and with Dirichlet boundary condition:
\begin{equation}
\begin{array}{cc}
-\nabla\left ( \kappa(x,\bZ) \nabla u(x,\bZ)\right)  = f(x), & x \in \Omega \subset \mathbb{R}^{2}, \\ 
u=g(x), & x \in \partial\Omega, 
\end{array}
\label{eq:elliptic}
\end{equation}
where $\kappa(x,\bZ)$ is a random field dependent on a random
parameter $\bZ=(Z_1,...,Z_{n_{z}}) \in \mathbb{R}^{n_{z}}$, $n_z\geq 1$, consisting of a set of independent continuous random variables characterizing the random coefficient of the governing equation. The solution  $u(x,\bZ)$ is a stochastic quantity, given by
\begin{equation}
u(x,\bZ):\overline{\Omega}\times \mathbb{R}^{n_{z}} \rightarrow \mathbb{R}^{n},
\end{equation}
where $n$ is the number of finite element nodes in $\Omega$.

For a fixed $\bZ$, the solution $u(x,\bZ)$ belongs to $H^1(\Omega)$, and for a fixed $x$ to $L_2(\vTheta)$.
There is an established theory about the existence and uniqueness of the solution to Eq.~\ref{eq:elliptic} under various assumptions on $\kappa$ and $f$; see, for example, 
\cite{babuska2004galerkin,Sarkis09,GITTELSON10,matthiesKeese05cmame,mugler2011elliptic}. In \cite{Sarkis09, GITTELSON10} it is shown that under additional
assumptions on the right-hand side $f$ and special choices of the test space the problem Eq.~\ref{eq:elliptic} is well-posed.
The case where the Lax-Milgram theorem is not applicable (e.g., upper and lower constants $\underline{\kappa}$, $\overline{\kappa}$ in
$0<\underline{\kappa}<\kappa<\overline{\kappa}<\infty$ do not exist)
is also considered in \cite{mugler2011elliptic}.
In \cite{Sarkis09} the authors analyze assumptions on $\kappa$ from \cite{babuska2004galerkin} to guarantee the uniqueness and the existence of the solution. Additionally, they offer a new method with weaker assumptions.
If the expansion of $\kappa$ is truncated, there is no guaranteed that the truncated series will stay strictly bounded from zero.
As a result, the existence of the approximate solution to Eq.~\ref{eq:elliptic} is questionable, unless precautions are taken as in \cite{matthiesKeese05cmame}.
The settings where the ellipticity condition is preserved are considered in \cite{Sarkis09}.

Further we assume that each continuous random variable $Z_i$ has a prior distribution
\begin{equation}
F_i(z_i)=P(Z_i \leq z_i)\in [0,1],
\end{equation}
where $P$ denotes probability and $\pp_i(z_i)=\frac{d F_i(z_i)}{d z_i}$ probability density function (pdf). The joint prior density function for $\bZ$ is $\pp_{\bZ}(z)=\prod _{i=1}^{n_{z}} \pp_i(z_i)$. For the sake of simplicity, we will skip the subscript ${}_{\bZ}$ and will write $\pp(z)$ for denoting the probability density function of the random variable $\bZ$. 
%

The elliptic boundary value problem in Eq.~\ref{eq:elliptic}, can represent, for instance, an incompressible single-phase porous media flow or, another example,  a steady state heat conduction through a composite material. In the single-phase flow, $u$ is the flow potential, and $\kappa$ the permeability of the porous medium. For heat conduction in composite materials, $u$ is the temperature, $-\kappa\nabla u$ the heat flow density, and $\kappa$ the thermal conductivity.

%

Iterative methods and preconditioners to solve the problem in Eq.~\ref{eq:elliptic} were developed in 
\cite{KressTobler11, KressTobler11HighDim, matthieszander-lowrank-2012, Sousedik2014truncated, ZanderDiss}.
In \cite{litvinenko-spde-2013} the authors assume that the solution has a low-rank \emph{canonical} (CP) tensor format and develop methods for the CP-formatted postprocessing.

Tensor ranks of the stochastic operator were analysed in \cite{matthies2012parametric,wahnert-stochgalerkin-2014}.
The proper generalized decomposition was applied for solving high dimensional stochastic problems in \cite{Nouy07,nouy-pgd-stoch-2010}.
In \cite{KhSch-Galerkin-SPDE-2011} authors employed newer tensor formats for the
approximation of coefficients and the solution of stochastic elliptic PDEs.
Other classical techniques to cope with high-dimensional problems are sparse grids
\cite{Griebel,Griebel_Bungartz, nobile2008sparse}
and (quasi) Monte Carlo methods \cite{graham-QMC-2011,scheichl-mlmc-further-2013,Schwab-MLQMC-2015}.
In \cite{dolgov2015polynomial, dolgov2014computation} authors approximate the polynomial chaos expansion (PCE ) of the random input coefficient $\kappa(x,\bZ)$ in the tensor train (TT) data format, and then solve the problem in that format. A low-rank tensor approximation of random fields, covariance matrices and set of snapshots is done in \cite{khoromskij2008data, LitvSampling13, litvinenko2013numerical}.

\subsection{Bayesian updating formula}
\label{S:bayes}

The inverse problem and 
propagation of uncertainty through a computational (forward)
model are strongly connected. Prior and posterior probabilities express our belief about possible values of the parameters $\kappa(x,\bZ)$ before and after observations.

Various ideas to speed up the Bayesian updating procedure were presented in 
\cite{Marzouk:2009, marzouk:2007:JCP, Najm:2008a, Berry:2012}. Surrogate based techniques were presented in \cite{rosic2013parameter,  matthies2012parametric, litvinenko2013inverse}; reduction of the stochastic dimension by using KLE and PCE expansions in \cite{rosic2012sampling, pajonk2012deterministic,rosic2013parameter};  a non-linear Kalman filter extension in \cite{matthies2016parameter, hermann2016inverse, litvinenko2013uncertainty}.

In \cite{Marzouk12OptMap}, the authors develop an approach to Bayesian inference that entirely avoids the Markov chain simulation by constructing a map that pushes forward the prior measure to the posterior measure. 
The work \cite{spantini2016goal} is devoted to optimal dimensionality reduction techniques for goal-oriented linear-Gaussian inverse problems, where the quantity of interest is a function of the inversion parameters.
A multiscale strategy for Bayesian inference using transport maps was introduced in \cite{parno2016multiscaleBayes}.

Further we assume that $\vTheta$
is a measure space with $\sigma$-algebra $\mathcal{A}$ and
with a probability measure $\mathbb{P}$, and that
$q: \vTheta \to \mathcal{Q}$ and $u: \vTheta \to \mathcal{U}$ are random variables (RVs).
Often, we are not able to observe
the entity $q\in \mathcal{Q}$ directly, we can only see
a `shadow' of it, formally given by a `measurement operator'
\begin{equation}  \label{eq:iI}
Y: \mathcal{Q} \times \mathcal{U} \ni (q,u) \mapsto  Y(q; u) \in \mathcal{Y},
\end{equation}
where $q(x,Z)=\log(\kappa(x,Z))$.
We assume that the space of possible measurements $\mathcal{Y}$ is a vector space,  
which frequently can be regarded as finite-dimensional,
as one can only observe a finite number of quantities.

The measurement operator $Y$ with values in $\mathcal{Y}$ produces
\begin{equation*}
y(Z) = Y(q(Z); u), \quad \text{where} \quad u=u(q(Z)). 
\end{equation*}
Examples of measurements are
a) $y(\bZ)=\int_{\omega}u(\bZ,x) dx$,  with a subdomain $\omega \subset \Omega$, and b) $u$ in a few points.
For a given $f$, the measurement $y$ is just a function of $q$. This function is usually not invertible since the measurement $y$ does not contain enough information.
In the Bayesian framework, the state of knowledge
is modeled in a probabilistic way. The parameter $q$ is uncertain and is modeled by a random variable.
The Bayesian setting allows updating/sharpening of
information about $q$ when the measurement is performed.

Usually the observation of the ``truth" $\hat{y} \in \mathbb{R}^{n_y}$
will deviate from what we expect to observe even if we know the right $q$
due to some model error
$\epsilon$. The measurement can be also polluted by some measurement error
$\vepsilon$.
Hence we observe $y= \hat{y}+ \epsilon + \vepsilon$, and would like to know what $q$ is. Let $\mathcal{S}:\mathbb{R}^{n_z}\rightarrow \mathbb{R}^{n_y}$ be the solution operator (for instance, the set $\{\Phi^g,\Phi^f\}$ or the inverse) of Eq.~\ref{eq:elliptic}. For the
sake of simplicity we will only consider one error term
\begin{equation}
y = \hat{y} + \vepsilon=\mathcal{S}(\bZ)+\vepsilon,\quad \text{where}\quad \vepsilon=(\vepsilon_1,...,\vepsilon_{n_y}) \in \mathbb{R}^{n_{y}}\;\; \text{includes all the errors}.
\end{equation}
Here $\vepsilon_1,...,\vepsilon_{n_y}$ are mutually independent random variables with probability density function $\pp(\vepsilon)=\prod_{i=1}^{n_{y}}\pp(\vepsilon_i)$. We also assume here that $\vepsilon$ and $\bZ$ are independent.

The mapping in \refeq{eq:iI} is usually not invertible, and hence the problem
is called ill-posed. By modeling our lack of knowledge
about $q$ in a Bayesian way \cite{Tarantola2004} 
with a $\mathcal{Q}$-valued random variable, the problem
becomes well-posed \cite{Stuart2010}.  But of course one is looking now at the
problem of finding a probability distribution that best fits the data;
and one also obtains a probability distribution of
$q$.  Here we focus on the use of Bayesian approach 
\cite{Goldstein2007}.

Bayes's theorem is commonly accepted as a consistent way to incorporate
new knowledge into a probabilistic description.
It may be formulated
as (\cite{Tarantola2004} Ch.\ 1.5)
\begin{equation}  \label{eq:iIIa}
 \pp(z|y) = \frac{\pp(y|z) }{\int_{\vTheta} \pp(y\vert z) \pp_z(z)\di z } \pp_z(z),
\end{equation}
where $\pp_z(z)$ is the pdf of $Z$, $\pp(y|z)$ is the likelihood as a function of $y$ for fixed prior $Z$ and $\pp(z|y)$ is the posterior pdf of $Z$ conditioned on the data $y$. We follow the notation from \cite{MarzoukXiu09}. Numerical approaches for computing a posterior
pdf were developed in \cite{marzouk:2007:JCP, Marzouk:2009, Stuart2010, rosic2013parameter}.
Assuming independence on the measurement noise $\varepsilon=(\varepsilon_1,\ldots, \varepsilon_{n_y})$, the likelihood function becomes
\begin{equation}
 L(z) := \pp(y|z)=\prod_{i=1}^{n_{y}}\pp_{\varepsilon_i}(y_i-\mathcal{S}_i(z)).
\end{equation}
Again, we see a formula, where the noisy measurement $y_i$ should be compared with the computed simulation $\mathcal{S}_i(z)$. And very often, the complete solution is not required.
\section{Hierarchical domain decomposition (HDD) method}
\label{sec:HDD}
The hierarchical domain decomposition (HDD) method \cite{MYPHD} combines the weak formulation, the finite element method (FEM), and the recursive domain decomposition method to obtain a fast and efficient algorithm for computing the partial inverse and a part of the solution (without computing the complete solution). 
This method was introduced by Hackbusch in 2002 and later on developed in \cite{MYPHD, HDDdoc, Hackbusch2012, FlorianDiss, HackHMBook}.

HDD computes the solution operators $\mathcal{F}_h$ and $\mathcal{G}_h$ in Eq.~\ref{eq:oper_represent2}, which after applying to the boundary condition and the right-hand side give us the solution.

\begin{figure}[h!]
\includegraphics[width=\textwidth]{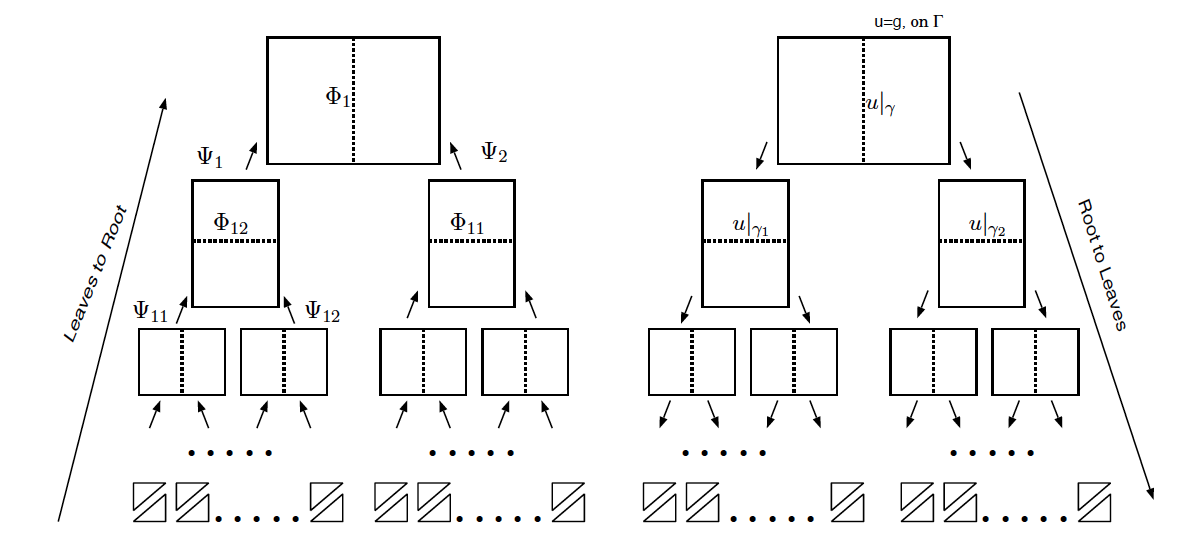}
\caption{HDD contains two algorithms: ``Leaves to Root'' (shown on the left) which computes mappings $\{\Psi_1, \Psi_2, \Psi_{11}, \Psi_{12},\ldots\}$ and $\{\Phi_1,\Phi_2,\Phi_{11},\Phi_{12},\ldots\}$ and ``Root to Leaves'' (on the right) which applies mappings $\{\Phi_{ij}\}$ to compute the solutions $u\vert_{\gamma_i}$ on the interfaces $\gamma_i$.}
\label{fig:HDD_th2}
\end{figure}

Below in this section we define the main components of the HDD method - the hierarchical domain decomposition tree (see Fig.~\ref{fig:HDD_th2}) in Section~\ref{sec:notationHDD}, the boundary-to-boundary mappings ($\Psi^g$) in Section~\ref{sec:Psi}, domain-to-boundary ($\Psi^f$) mappings in Section~\ref{sec:Phi}, boundary-to-interface ($\Phi^g$) and domain-to-interface ($\Phi^f$) mappings which are essential for the definition of the HDD method in Section~\ref{sec:Phi}.

For a fixed parameter $Z$, Eq.~\ref{eq:elliptic} can be written as follow: 
\begin{align}
\label{eq:elliptic_det}
-\nabla \left( \kappa(\x) \nabla u(\x)\right) &= f(\x), \quad \x \in \Omega \subset \mathbb{R}^{2},\\ 
u&=g(\x), \quad \x \in \partial\Omega, \nonumber
\end{align}
where $\x=(x_1,x_2)\in \Omega$.

The HDD method computes two discrete hierarchical solution operators $\F_h$ and $\G_h$ such that:
\begin{equation}
\label{eq:oper_represent2}
u_h=\F_hf_h+\G_hg_h,  
\end{equation}
where $u_h=u_h(f_h,g_h)$ is the FE solution of \ref{eq:elliptic_det}, $f_h$ the discretized right-hand side, and $g_h$ the Dirichlet boundary data. To decrease the computing time and the storage cost, both operators $\F_h$ and $\G_h$ are approximated by $\H$-matrices.

\begin{figure}[htbp]
\centerline{\includegraphics[width=5.5in, height=1.9in]{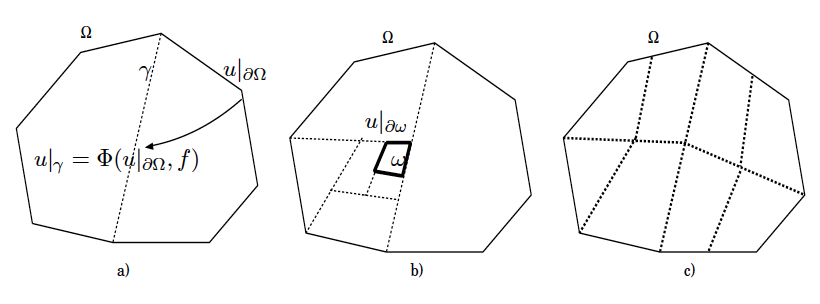}}
\caption{(a) The solution  $u\vert_{\gamma}$ on the interface $\gamma$ can be computed with the auxiliary operator $\Phi$, by applying it to the right hand side $f$ and to the boundary condition $u\vert_{\partial\Omega}$; (b) HDD can compute the solution in a subdomain $\omega \subset \Omega$; (c) HDD method can compute the solution on a coarse mesh (shown by dotted lines).}
\label{fig:HDDidea}
\end{figure} 


\textbf{Three examples of possible problem setups, shown in Fig.~\ref{fig:HDDidea}}, are the following:
\begin{enumerate}
\item Suppose the solution on the boundary $\partial\Omega$ (Fig. \ref{fig:HDDidea} (a)) is given.
One is interested in the fast numerical approach which computes the solution $u\vert_{\gamma}$ on the interface $\gamma$. The solution $u\vert_{\gamma}$ depends on the right-hand side and $u\vert_{\partial\Omega}=g$, i.e.  $u\vert_{\gamma}=\Phi(u\vert_{\partial\Omega}, f)$ with some mapping $\Phi$;
\item
Only the solution in a small subdomain $\omega \subset \Omega$ is of interest (Fig. \ref{fig:HDDidea} (b)). To solve the problem in a domain $\omega$ the boundary values on $\partial\omega$ are required. How to compute them efficiently from the global boundary data $\partial\Omega$ and the given right-hand side? 
\item
The third possible problem setup is as follows. The solution on the interface or on a very coarse mesh (see Fig. \ref{fig:HDDidea} (c)) is required. How can this solution be computed effectively without neglecting small scale features?
\end{enumerate}
Other properties of the HDD method are the following. The HDD allows one to compute $u_h(f_h,g_h)$ for $f_h$ given in a smaller space $V_H \subset V_h$. This could be useful, for instance, in multi-scale settings.
The HDD provides the possibility to compute $u_h$ restricted to a coarser grid with reduced  computational effort. 
The HDD shows big advantages in complexity for problems with multiple right-hand sides and multiple Dirichlet data. In this case both operators $\F_h$ and $\G_h$ are computed only once and then applied multiple times to $f_h$ and $g_h$.    
Due to the binary tree structure the HDD is an easily parallelizable method.
If the problem contains repeated patterns (for instance, so-called \textit{cells} in a multi-scale framework) then the computational resources can be reduced drastically.
\subsection{Notation}
\label{sec:notationHDD}
Let $\Th$ be a triangulation of the spatial domain $\Omega$. 
After hierarchical decomposition of $\Omega$ (cf. \cite{NestedDissection}), obtain the \textit{hierarchical domain decomposition tree} $T_{\mathcal{T}_h}$  (see Fig.~\ref{fig:HDD_th2}) with the following properties:
\begin{itemize}
\item[$\bullet$] $\Omega $ is the root of the tree,

\item[$\bullet$] $T_{\mathcal{T}_h} $ is a binary tree,

\item[$\bullet$] If $\omega \in T_{\mathcal{T}_h} $ has two sons 
$\omega_1 ,\omega_2
\in T_{\mathcal{T}_h}$, then \\
$\omega = \omega_1 \cup \omega_2 $ and $%
\omega_1 ,\omega_2 $ have no interior point in common,
\item [$\bullet$]
$\omega \in T_{\mathcal{T}_h} $ is a leaf, if and only if $\omega \in 
\mathcal{T}_h$.
\end{itemize}
The construction of $T_{\Th}$ is straight-forward by dividing $\Omega$ recursively into subdomains. For practical purposes, the subdomains $\omega_1 $, $\omega_2$ must both be of size $\approx \vert \omega \vert /2$ and the internal boundary
\begin{equation}
\label{eq:gamma_def}
\gamma_{\omega}:=\partial\omega_1 \backslash \partial\omega=\partial\omega_2 \backslash \partial\omega
\end{equation}
must not be too large (see Fig. \ref{fig:devide} (left)).

\begin{figure}[htb]
\begin{center}
{\includegraphics[width=0.29\textwidth]{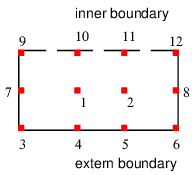}}
\end{center}
\caption{Domain $\omega_1\in T_{\Th}$ with $I({\omega_1})=\{1,...,12\}$, $I({\partial \omega_1})=\{3,4,5,6,7,8,9,10,11,12\}$, $I({\gamma_{\omega_1}})=\{1,2\}$ to be eliminated via the Schur complement, $I({\Gamma_{\omega_1}})=\{9,7,3,4,5,6,8,12\}$. On the next level, when $\omega_1$ will be coupled with $\omega_2$, the points $I({\gamma_\omega})=\{10,11\}$ will be eliminated.}
\label{fig:ex_points}
\end{figure}

Let $I:=I(\overline{\Omega})$ and $x_i$, $i\in I$, be the {set of all nodal points} in $\overline{\Omega}$ (including nodal points on the boundary). We define $I(\omega)$ as a subset of $I$ with $x_i\in \omega=\overline{\omega}$. Similarly, we define $I(\mathop \omega \limits^{\circ})$, $I(\Gamma_{\omega})$, $I(\gamma_{\omega})$, where $\Gamma_{\omega}:=\partial\omega$, $\mathop \omega \limits^{\circ}=\omega\backslash \partial\omega$, for the interior, for the external boundary and for the interface.

Computing the discrete solution $u_h$, Eq.~\ref{eq:elliptic}, in $\Omega$ is equivalent to the computation of $u_h$ on all $\gamma_{\omega}$, $\omega \in T_{\mathcal{T}_h}$, since $\displaystyle{I(\Omega)=\cup_{\omega \in T_{\mathcal{T}_h}}I(\gamma_{\omega})}$. These computations are performed by
using the linear mappings $\Phi_{\omega}^f$, $\Phi_{\omega}^g$ defined for all nodes ${\omega \in T_{\mathcal{T}_h}}$.
\begin{notation}
\label{not:g_f}
Let $g_{\omega}:=u\vert_{I(\partial\omega)}$ be the local Dirichlet data and $f_{\omega}:=f\vert_{I(\omega)}$ be the local right-hand side.
\end{notation}
\begin{defi}
The mapping $\Phi^g_{\omega}:\RR^{I(\partial\omega)}\rightarrow \RR^{I(\gamma_{\omega})}$ maps the boundary data defined on $\partial\omega$ to the data defined on the interface $\gamma_{\omega}$. $\Phi^f_{\omega}:\RR^{I(\omega)}\rightarrow \RR^{I(\gamma_{\omega})}$ maps the right-hand side data defined on $\omega$ to the data defined on $\gamma_{\omega}$.\\
\end{defi}
The final aim is to compute the solution $u_h$ along $\gamma_{\omega}$ in the form 
$u_h\vert_{\gamma_{\omega}}=\Phi^f_{\omega}f_{\omega}+\Phi^g_{\omega}g_{\omega}$, $\omega \in T_{\Th}$. 
For this purpose HDD builds the mappings $\Phi_{\omega}:=(\Phi^g_{\omega}, \Phi^f_{\omega})$, for all $\omega \in T_{\mathcal{T}_h}$. For computing the mapping $\Phi_{\omega}$, $\omega \in T_{\Th}$, we first need to compute the auxiliary mapping $\Psi_{\omega}:=(\Psi^g_{\omega},\Psi^f_{\omega})$ which will be defined later.

Thus, the HDD method consists of two steps: the first step is the construction of
the mappings $\Phi^g_{\omega}$ and $\Phi^f_{\omega}$ for all $\omega \in T_{\mathcal{T}_h}$. The second step is the recursive computation of the solution $u_h$. In the second step HDD applies the mappings $\Phi^g_{\omega}$ and $\Phi^f_{\omega}$ to 
the local Dirichlet data $g_{\omega}$ and to the local right-hand side $f_{\omega}$.
%
\begin{notation}
Let $\omega \in T_{\mathcal{T}_h}$ and
\begin{equation}
\label{eq:d_definit}
d_{\omega}: = \left( {\left( {f_i} \right)_{i \in I(\omega)} ,\left( {g_i }
\right)_{i \in I(\partial_{\omega} )} } \right)=(f_{\omega},g_{\omega})
\end{equation}
be a composed vector consisting of the right-hand side from Eq.~\ref{eq:elliptic} restricted to $\omega$ and the Dirichlet boundary values $g_{\omega}=u_h\vert_{\partial\omega}$ (see also Notation \ref{not:g_f}).
\end{notation}
Note that $g_{\omega}$ coincides with the global Dirichlet data in Eq.~\ref{eq:elliptic} only when $\omega=\Omega$. For all other $\omega \in T_{\Th}$ we compute $g_{\omega}$ in (Eq.~\ref{eq:d_definit}) by the algorithm ``Root to Leaves'' (see Section \ref{sec:root2leaves}).

Assuming that the elliptic boundary
value problem, Eq.~\ref{eq:elliptic_det}, restricted to $\omega$ is solvable, we can define the local FE solution by solving the following discrete problem in the variational form \cite{HackEDE}: 
\begin{equation} \label{eq:bilform}
\left\{ 
\begin{array}{ll}
a_\omega(U_{\omega},b_j) = \left({f_{\omega},b_j} \right)_{L^2(\omega)}, & \forall \;
j \in I(\mathop \omega \limits^{\circ}), \\ 
U_{\omega} (\x_j ) = g_j, & \forall \; j \in I(\partial \omega).%
\end{array}
\right.
\end{equation}
Here, $b_j $ is the $P^1$-Lagrange basis function at $\x_j$ 
and $a_\omega ( \cdot , \cdot )$ is the bilinear form (see Eq.~\ref{eq:elliptic}) with integration
restricted to $\omega$ and $(f_{\omega},b_j)=\int\limits_\omega {f_{\omega}\,b_j\,d\x 
}$. \\
Let $U_{\omega}\in V_h$ be the solution of (Eq.~\ref{eq:bilform}) in $\omega$. The solution $U_{\omega}$ depends on the Dirichlet data on $\partial\omega$ and the right-hand side in $\omega$. Dividing problem (Eq.~\ref{eq:bilform}) into two subproblems (Eq.~\ref{eq:bilform1}) and (Eq.~\ref{eq:bilform2}), we obtain $U_{\omega}=U_{\omega}^f+U_{\omega}^g$, where $U_{\omega}^f$ is the solution of
\begin{equation}  \label{eq:bilform1}
\left\{ 
\begin{array}{ll}
a_\omega(U_{\omega}^f,b_j) = \left({f_{\omega},b_j} \right)_{L^2(\omega)}, & \forall \;
j \in I(\mathop \omega \limits^{\circ}), \\ 
U_{\omega}^f (\x_j ) = 0, & \forall \; j \in I(\partial \omega)%
\end{array}
\right.
\end{equation}
and $U_{\omega}^g$ is the solution of
\begin{equation}  \label{eq:bilform2}
\left\{ 
\begin{array}{ll}
a_\omega(U_{\omega}^g,b_j) = 0, & \forall \;
j \in I(\mathop \omega \limits^{\circ}), \\ 
U_{\omega}^g (\x_j ) = g_j, & \forall \; j \in I(\partial \omega).%
\end{array}
\right.
\end{equation}
If $\omega=\Omega$ then (Eq.~\ref{eq:bilform}) is equivalent to the initial problem Eq.~\ref{eq:elliptic_det} in the weak formulation.
\subsection{Mapping $\Phi_{\omega}=(\Phi_{\omega}^g,\Phi_{\omega}^f)$}
\label{sec:Phi}
In this section we define mappings $\Phi_{\omega}$, $\Phi_{\omega}^g$, $\Phi_{\omega}^f$.
We consider $\omega \in T_{\mathcal{T}_h} $ with two sons $\omega_1 ,\omega_2
$. Considering once more the data $d_\omega$ from (Eq.~\ref{eq:d_definit}), $U^f_{\omega}$ from (Eq.~\ref{eq:bilform1}) and $U^g_{\omega}$ from (Eq.~\ref{eq:bilform2}), we define $\Phi^f_\omega (f_\omega)$ and $\Phi^g_\omega (g_\omega)$ by
\begin{equation}
\label{eq:fi_f}
\left( {\Phi^f_\omega(f_{\omega})} \right)_i : = U^f_{\omega} (\x_i )\, \quad \forall i \in I(\gamma_\omega )
\end{equation}
and
\begin{equation}
\label{eq:fi_g}
\left( {\Phi^g_\omega(g_{\omega})} \right)_i : = U^g_{\omega} (\x_i )\, \quad \forall i \in I(\gamma_\omega ).
\end{equation}
Since $U_{\omega}=U^f_{\omega}+U^g_{\omega}$, we obtain
\begin{equation}
\label{eq:fi}
\left( {\Phi_\omega(d_{\omega})} \right)_i := \Phi^g_\omega(g_{\omega})+\Phi^f_\omega(f_{\omega})=U^f_{\omega} (\x_i )+U^g_{\omega} (\x_i )=U_{\omega} (\x_i )
\end{equation}
for all $i \in I(\gamma_\omega )$.\\
\noindent Hence, $\Phi_\omega(d_{\omega})$ is the trace of $U_{\omega}$ on $\gamma _\omega $.
Definition in (Eq.~\ref{eq:fi}) says that if the data $d_{\omega}$ are given then $\Phi_{\omega}$ computes the solution of (Eq.~\ref{eq:bilform}). Indeed, $\Phi_{\omega}d_{\omega}=\Phi^g g_{\omega} + \Phi^f f_{\omega}$. Note that the solution $u_h$ of the initial global problem coincide with $U_{\omega}$ in $\omega$, i.e., $u_h\vert_{\omega}=U_{\omega}$. 
\subsection{Mapping $\Psi_{\omega}=(\Psi_{\omega}^g, \Psi_{\omega}^f) $} 
\label{sec:Psi}
In this section we define mappings $\Psi_{\omega}$, $\Psi_{\omega}^g$, $\Psi_{\omega}^f $.

First, we define the mapping $\Psi^f_{\omega}$ from (Eq.~\ref{eq:bilform1}) as
\begin{equation}  \label{eq:psi_f}
\begin{array}{l}
\left( {\Psi_{\omega}^f (d_\omega)} \right)_{i \in I(\partial\omega)} := a_\omega
(U_{\omega}^f ,b_i ) - \left( {f_{\omega},b_i } \right)_{L^2(\omega)},%
\end{array}%
\end{equation}
where $U_{\omega}^f \in V_h$, $U_{\omega}^f\vert_{\partial\omega}=0$ and
\begin{equation*}
a(U_{\omega}^f, b_i)-(f,b_i)=0, \quad \text{for }\forall i\in I(\wo).
\end{equation*}
Second, we define the mapping $\Psi^g_{\omega}$ from (Eq.~\ref{eq:bilform2}) by setting
\begin{equation}  \label{eq:psi_g}
\begin{array}{l}
\left( {\Psi_{\omega}^g (d_\omega)} \right)_{i \in I(\partial
\omega)} := a_\omega(U_{\omega}^g ,b_i ) - \left( {f_{\omega},b_i } \right)_{L^2(\omega)}=
a_\omega(U_{\omega}^g ,b_i )-0=a_\omega(U_{\omega}^g ,b_i ),
\end{array}%
\end{equation}
where $U_{\omega}^g \in V_h$ and $\left( {\Psi_{\omega}^g (d_\omega)} \right)_i=0$ for $\forall i\in I(\wo)$.

The linear mapping $\Psi_{\omega}$, which maps the
data $d_{\omega}$ given by (Eq.~\ref{eq:d_definit}) to the boundary
data on $\partial \omega$, is given in the component form as
\begin{equation}  \label{eq:psi}
\begin{array}{l}
\Psi_\omega (d_{\omega}) = \left( {\Psi_\omega (d_\omega)} \right) _{i \in I(\partial
\omega)} := a_\omega
(U_{\omega} ,b_i ) - \left( {f_{\omega},b_i } \right)_{L^2(\omega)}.%
\end{array}%
\end{equation}
By definition $\Psi_\omega $ is linear in $(f_{\omega},g_{\omega})$ and can be written as $%
\Psi_\omega (d_{\omega}) = \Psi_{\omega}^f f_{\omega} + \Psi_{\omega}^g g_{\omega}$.
Here $U_{\omega}$ is the solution of the local problem (Eq.~\ref{eq:bilform}) and it coincides with the global solution on $I(\omega)$. 
%
%
\subsection{$\Phi_{\omega}$ and $\Psi_{\omega}$ in terms of the Schur complement matrix}
\label{sec:mapsSchur}

Let the linear system $A\u=F\c$ for $\omega \in T_{\Th}$ be given. In Sections \ref{sec:start_recursion} and \ref{sec:build_father} we explain how to obtain the matrices $A$ and $F$. $A$ is the stiffness matrix for the domain $\overline{\omega}$ after elimination of the unknowns corresponding to $I(\mathop {\omega} \limits^{\circ} \setminus \gamma_{\omega})$. The matrix $F$ comes from the applied numerical integration rule \cite{MYPHD}.

We will write for simplicity $\gamma$ instead of $\gamma_{\omega}$.
Thus, $A:\RR^{I(\partial\omega\cup \gamma)}\rightarrow \RR^{I(\partial\omega\cup \gamma)}$, $\u\in\RR^{I(\partial\omega\cup \gamma)}$, 
$F:\RR^{I(\omega)}\rightarrow \RR^{I(\partial\omega\cup \gamma)}$ and 
$\c\in \RR^{I(\omega)}$. Decomposing the unknown vector $\u$ into 
two components $\u_1\in\RR^{I(\partial\omega)}$ and 
$\u_2\in\RR^{I(\gamma)}$, obtain
\begin{displaymath}
\u=\left(\begin{array}{l}
\u_1 \\ 
\u_2
\end{array}\right).
\end{displaymath}
The component $\u_1$ corresponds to the boundary $\partial\omega$ and the component $\u_2$ to the interface $\gamma$.
Then the equation $A\u=F\c$ becomes
\begin{equation}
\label{eq:block_equation1}
\left(\begin{array}{cc}
A_{11} & A_{12} \\ 
A_{21} & A_{22}
\end{array}\right)
\left(\begin{array}{c}
\u_1 \\ 
\u_2
\end{array}\right)
=
\left(\begin{array}{c}
F_1 \\ 
F_2
\end{array}\right)\c,
\end{equation}
where
\begin{equation*}
A_{11}:\RR^{I(\partial \omega)}\rightarrow \RR^{I(\partial \omega)}, \quad
A_{12}:\RR^{I(\gamma)}\rightarrow \RR^{I(\partial \omega)},
\end{equation*}
\begin{equation*}
A_{21}:\RR^{I(\partial\omega)}\rightarrow \RR^{I(\gamma)},\quad
A_{22}:\RR^{I(\gamma)}\rightarrow \RR^{I(\gamma)},
\end{equation*}
\begin{equation*}
F_{1}:\RR^{I(\omega)}\rightarrow \RR^{I(\partial \omega)},\quad
F_{2}:\RR^{I(\omega)}\rightarrow \RR^{I(\gamma)}.
\end{equation*}

The elimination of the internal points is done as it is shown in (Eq.~\ref{eq:SchurMinus1}) below
\begin{equation}
\label{eq:SchurMinus1}
\left(\begin{array}{cc}
A_{11}-A_{12}A_{22}^{-1}A_{21} & 0 \\ 
A_{21} & A_{22}
\end{array}\right)
\left(\begin{array}{c}
\u_1 \\ 
\u_2
\end{array}\right)
=
\left(\begin{array}{c}
F_1-A_{12}A_{22}^{-1}F_2 \\ 
F_2
\end{array}\right)\c.
\end{equation}
We rewrite the last system as two equations
\begin{equation}
\label{eq:SchurC1}
\begin{array}{c}
\tilde{A}\u_1:=(A_{11}-A_{12}A_{22}^{-1}A_{21})\u_1=(F_1-A_{12}A_{22}^{-1}F_2)\c,
 \\ 
\u_2=A_{22}^{-1}F_2\c-A_{22}^{-1}A_{21}\u_1.
\end{array} 
\end{equation}
The explicit expressions for the mappings $\Psi_{\omega}$ and $\Phi_{\omega}$ follow from (Eq.~\ref{eq:SchurC1}): 
\begin{equation}
\Psi^g_{\omega}:=A_{11}-A_{12}A_{22}^{-1}A_{21},\quad
\Psi^f_{\omega}:=F_1 - A_{12}A_{22}^{-1}F_2,
\end{equation}
\begin{equation}
\Phi^g_{\omega}:=-A_{22}^{-1}A_{21},\quad
\Phi^f_{\omega}:=A_{22}^{-1}F_2.
\end{equation}
Thus, $\u_2=\Phi^f_{\omega}(f_{\omega})+\Phi^g_{\omega}(g_{\omega})$, with the rhs $f_{\omega}=\c$, and local b.c. $g_{\omega}=\u_1$.
\section{Construction Process}
\label{sec:HDDConstrProcess}
In this section we explain the recursive construction of mappings 
$\Psi^g_{\omega}$, $\Psi^f_{\omega}$,
$\Phi^g_{\omega}$ and $\Phi^f_{\omega}$.

\subsection{Initialisation of the recursion}
\label{sec:start_recursion}
This section explains how to compute mapping $\Psi^f_{\omega}$ for the leaves of $T_{\Th}$ and how it is connected with the quadrature rule.

Our purpose is to get for each triangle $\omega \in \Th$, the system of linear equations
\begin{equation}
\label{eq:AuFc}
A\cdot \u=\tilde{\c}:=F\cdot \c,
\end{equation}
where $A$ is the stiffness matrix, $\c$ the discrete values of the right-hand side in the nodes of $\omega$ and $F$ will be defined later. The matrix coefficients $A_{ij}$ are computed by the formula
\begin{equation}  \label{eq:coefA}
A_{ij} = \int\limits_{\omega}{\kappa(\x)\langle \nabla b_i(\x) \cdot \nabla b_j(\x) \rangle d\x},
\end{equation}
where $b_i(\x)$ is a piecewise linear basis function \cite{HackEDE}.
For $\omega \in \Th$, $F\in \RR^{3\times 3}$ comes from the discrete integration and
the matrix coefficients $F_{ij}$ are computed using (Eq.~\ref{eq:coefF}).
The components of $\tilde{\c}$ can be computed as follows:
\begin{equation}  
\label{eq:coefc}
\tilde{c}_{i} = \int\limits_{\omega}{fb_id\x}\approx \frac{f(\x_1)b_i(\x_1)+f(\x_2)b_i(\x_2)+f(\x_3)b_i(\x_3)}{3}\cdot \vert \omega \vert,
\end{equation}
where $\x_i$, $i \in \{1,2,3\}$, are three vertices of the triangle $\omega \in T_{\Th}$, $b_i(\x_j)=1$ if $i=j$ and $b_i(\x_j)=0$ otherwise. Rewrite (Eq.~\ref{eq:coefc}) in matrix form:
\begin{equation}  
\label{eq:coefFv}
\tilde{\c}=
\left(\begin{array}{c}
\tilde{c}_1\\
\tilde{c}_2\\
\tilde{c}_3
\end{array}
\right)
\approx\frac{1}{3}
\left(\begin{array}{c c c}
b_1(\x_1) & b_1(\x_2)& b_1(\x_3)\\
b_2(\x_1) & b_2(\x_2)& b_2(\x_3)\\
b_3(\x_1) & b_3(\x_2)& b_3(\x_3)
\end{array}
\right)
\left(\begin{array}{c}
f(\x_1)\\
f(\x_2)\\
f(\x_3)
\end{array}
\right),
\end{equation} where $f(\x_i)$, $i=1,2,3$, are the values of the right-hand side $f$ in the vertices of $\omega$. Then, for piecewise linear basis functions obtain
\begin{equation}
\label{eq:coefF}
F:=
\frac{1}{3}
\left(\begin{array}{c c c}
b_1(\x_1) & b_1(\x_2)& b_1(\x_3)\\
b_2(\x_1) & b_2(\x_2)& b_2(\x_3)\\
b_3(\x_1) & b_3(\x_2)& b_3(\x_3)
\end{array}
\right)
=
\frac{1}{3}
\left(\begin{array}{c c c}
1 & 0& 0\\
0 & 1& 0\\
0 & 0& 1
\end{array}
\right)
\text{   and   }
c:=\left(\begin{array}{c}
f(\x_1)\\
f(\x_2)\\
f(\x_3)
\end{array}
\right).
\end{equation}
Thus, $\Psi^g_{\omega}$ corresponds to the matrix $A\in \RR^{3\times 3}$ and $\Psi^f_{\omega}$ to $F \in \RR^{3 \times 3}$. 
\subsection{Recursion}
\begin{figure}[htb]
\begin{center}
{\includegraphics[width=0.99\textwidth]{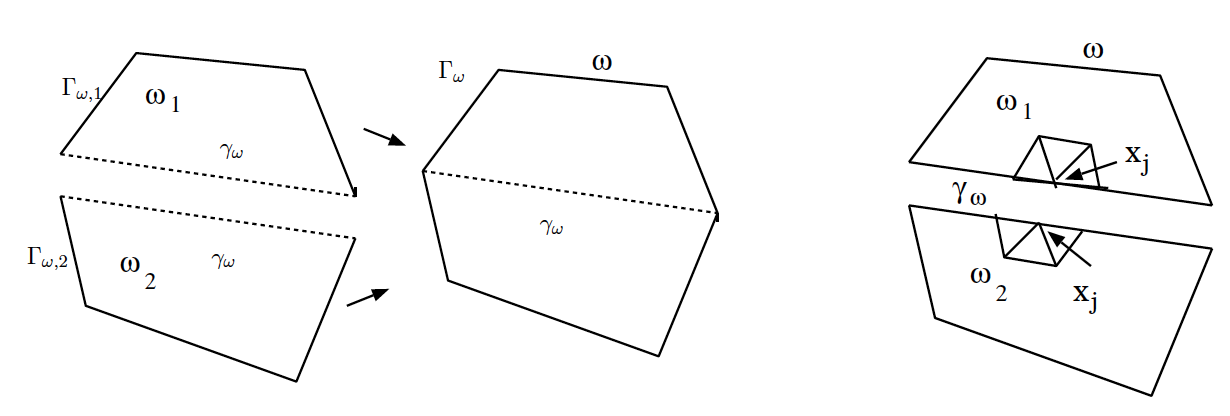}}
\end{center}
\caption{(\textbf{left}) Domain $\omega$ and its two sons $\omega_1$ and $\omega_2$. Here $\gamma_{\omega}$ is the internal boundary and $\Gamma_{\omega,i}$, $i=1,2$, parts of the external boundaries, see (Eq.~\ref{eq:Gamma12}). (\textbf{right}) The support of basis function $b_j$, $x_j\in \omega_1$ and $x_j\in \omega_2$.}
\label{fig:devide}
\end{figure}

This section explains how to build $\Psi_{\omega}$ from $\Psi_{\omega_1}$ and $\Psi_{\omega_2}$, with $\omega \in T_{\mathcal{T}_h} $ and $\omega_1$, $\omega_2$ be two sons of $\omega$.
The coefficients of $\Psi_\omega$ can be computed by (Eq.~\ref{eq:psi}). 
The external boundary $\Gamma_{\omega}$ of $\omega$ splits into (see Fig. \ref{fig:devide} (left))
\begin{equation}
\label{eq:Gamma12}
\Gamma_{\omega,1}:=\partial\omega \cap \omega_1, \quad \Gamma_{\omega,2}:=\partial\omega \cap \omega_2.
\end{equation}
For simplicity of further notations, we will write $\gamma$ instead of $\gamma_{\omega}$. %
\begin{notation}
Recall that $I(\partial\omega_i)=I(\Gamma_{\omega,i})\cup I(\gamma)$. We denote the restriction of $\Psi_{\omega_i}:\RR^{I(\partial\omega_i)}\rightarrow \RR^{I(\partial\omega_i)}$ to
$I(\gamma)$ by ${\ }^\gamma \Psi_{\omega}:=(\Psi_{\omega})\vert_{i \in I(\gamma)}$. 
\end{notation}
Suppose that by induction, the mappings $\Psi_{\omega_1 }$, $\Psi_{\omega_2 } $ are known
for the sons $\omega_1$, $\omega_2 $. Now, we explain how to construct $\Psi_\omega $ and $\Phi_\omega$. 

\begin{lemma}
Let the data $d_1 = d_{\omega_1}$, $d_2 = d_{\omega_2}$ be given by (Eq.~\ref%
{eq:d_definit}). Data $d_1$ and $d_2$ coincide along $\gamma _\omega $, i.e.,\\
$\bullet$ (consistency conditions for the boundary) 
\begin{equation}  \label{eq:dirich}
g_{1,i} = g_{2,i}\, \quad \forall i \in I(\omega_1 ) \cap I(\omega_2), 
\end{equation}
$\bullet$ (consistency conditions for the right-hand side) 
\begin{equation}  \label{eq:contin}
f_{1,i} = f_{2,i}\, \quad \forall i \in I(\omega_1 ) \cap I(\omega_2 ). 
\end{equation}
If the local FE solutions $u_{h,1} $ and $u_{h,2} $ of the problem (\ref%
{eq:bilform}) for the data $d_1 ,d_2 $ satisfy the additional equation  
\begin{equation}  
\label{eq:neumann}
{\ }^\gamma \Psi_{\omega_1 } (d_1 ) + {\ }^\gamma \Psi_{\omega_2} (d_2 ) = 0,
\end{equation}
\noindent then the composed solution $u_h $ defined by assembling
\begin{equation}  \label{eq:uh}
u_h (\x_i ) = \left\{ {{\begin{array}{*{20}c} {u_{h,1} (\x_i ) \quad
\text{for} \quad i \in I(\omega_1 ),} \hfill \\ {u_{h,2} (\x_i ) \quad
\text{for} \quad i \in I(\omega_2 )} \hfill \\ \end{array} }} \right.
\end{equation}
\noindent satisfies (Eq.~\ref{eq:bilform}) for the data $d_\omega = (f,g)$ where 
\begin{equation}  \label{eq:f}
f_i = \left\{ {{\begin{array}{*{20}c} {f_{1,i} \quad \text{for} \quad i \in
I(\omega_1 ),} \hfill \\ {f_{2,i} \quad \text{for} \quad i \in I(\omega_2
),} \hfill \\ \end{array} }} \right.
\end{equation}
\begin{equation}  \label{eq:g}
g_i = \left\{ {{\begin{array}{*{20}c} {g_{1,i} \quad \text{for} \quad i \in
I(\Gamma _{\omega,1} ),} \hfill \\ {g_{2,i} \quad \text{for} \quad i \in
I(\Gamma _{\omega,2} ).} \hfill \\ \end{array} }} \right.
\end{equation}
\end{lemma}
\beginproof
Note that the index sets in (Eq.~\ref{eq:uh})-(Eq.~\ref{eq:g})
overlap. Let $\omega_1 \in T_{\Th}$, $f_{1,i}=f_i$, $i\in I(\omega_1)$, and $g_{1,i}=g_i$, $i\in I(\partial\omega_1)$. Then the existence of the unique solutions of (Eq.~\ref{eq:bilform}) 
gives 
$u_{h,1}(\x_i)=u_{h}(\x_i)$, $\forall i \in I(\mathop {\omega_1} \limits^{\circ})$.\\
In a similar manner we get $u_{h,2}(\x_i)=u_{h}(\x_i)$ , $\forall i \in I(\mathop {\omega_2} \limits^{\circ})$.
Equation (Eq.~\ref{eq:psi}) gives 
\begin{equation} 
\begin{array}{l}
\left( {\ }^\gamma {\Psi_{\omega_1} (d_1)} \right)_{i \in I(\gamma)} = a_{\omega_1}
(u_h ,b_i ) - \left( {f_{\omega_1},b_i } \right)_{L^2(\omega_1)}
\end{array}
\end{equation}
and
\begin{equation} 
\begin{array}{l}
\left( {\ }^\gamma {\Psi_{\omega_2} (d_2)} \right)_{i \in I(\gamma)} = a_{\omega_2}
(u_h ,b_i ) - \left( {f_{\omega_2},b_i } \right)_{L^2(\omega_2)}.
\end{array}
\end{equation}
The sum of the two last equations (see Figure \ref{fig:devide} (right)) and (Eq.~\ref{eq:neumann}) give 
\begin{equation} 
\begin{array}{l}
0={\ }^\gamma {\Psi_\omega (d_{\omega})}_{i \in I(\gamma)} = a_\omega
(u_h ,b_i ) - ( {f_{\omega},b_i })_{L^2(\omega)}.
\end{array}
\end{equation}
We see that $u_h$ satisfies (Eq.~\ref{eq:bilform}). 
\endproof
Note that
\begin{equation*}
u_{h,1} (\x_i) = g_{1,i} = g_{2,i} = u_{h,2}(\x_i) \quad \text{holds for} \quad i \in I(\omega_1 ) \cap I(\omega_2 ).
\end{equation*}
Next, we use the decomposition of the data $d_1$ into the components 
\begin{equation}
\label{eq:d_components}
d_1 = (f_1, g_{1,\Gamma } ,g_{1,\gamma } ),
\end{equation}
where 
\begin{equation}
\label{eq:g_components}
g_{1,\Gamma } := (g_1 )_{i \in I(\Gamma _{\omega,1} )}, \quad g_{1,\gamma } :=
(g_1 )_{i \in I(\gamma)}
\end{equation}
and similarly for $d_2 = (f_2 ,g_{2,\Gamma } ,g_{2,\gamma } )$. \\
The decomposition $g \in \RR^{I(\partial\omega_j)}$ into $g_{j,\Gamma} \in \RR^{I(\Gamma_{\omega,j})}$ and $g_{j,\gamma} \in \RR^{I(\gamma)}$ implies the decomposition of $\Psi^g_{\omega_j}:\RR^{ I(\partial\omega_j)}\rightarrow \RR^{ I(\partial\omega_j)}$ into $\Psi^{\Gamma}_{\omega_j}:\RR^{ I(\Gamma_{\omega,j})}\rightarrow \RR^{ I(\partial\omega_j)}$ and $\Psi^{\gamma}_{\omega_j}:\RR^{ I(\gamma)}\rightarrow \RR^{ I(\partial\omega_j)}$, $j=1,2$. 
Thus, $\Psi^{g}_{\omega_1}g_{\omega_1}=\Psi^{\Gamma}_{\omega_1}g_{1,\Gamma}+\Psi^{\gamma}_{\omega_1}g_{1,\gamma}$ and $\Psi^{g}_{\omega_2}g_{\omega_2}=\Psi^{\Gamma}_{\omega_2}g_{2,\Gamma}+\Psi^{\gamma}_{\omega_2}g_{2,\gamma}$.\\
The maps $\Psi_{\omega_1}$, $\Psi_{\omega_2}$ become 
\begin{equation}
\label{eq:Psi1_zerleg}
\Psi_{\omega_1 } d_1 = \Psi_{\omega_1 }^f f_1 + \Psi_{\omega_1 }^\Gamma
g_{1,\Gamma } + \Psi_{\omega_1 }^\gamma g_{1,\gamma },
\end{equation}
\begin{equation} 
\label{eq:Psi2_zerleg}
\Psi_{\omega_2 } d_2 = \Psi_{\omega_2 }^f f_2 + \Psi_{\omega_2 }^\Gamma
g_{2,\Gamma } + \Psi_{\omega_2 }^\gamma g_{2,\gamma }.
\end{equation}
\begin{defi}
We will denote the restriction of $\Psi^{\gamma}_{\omega_j}:\RR^{ I(\gamma)}\rightarrow \RR^{ I(\partial\omega_j)}$ to $I(\gamma)$ by
\begin{equation*}
{\ }^\gamma\Psi^{\gamma}_{\omega_j}:\RR^{ I(\gamma)}\rightarrow \RR^{ I(\gamma)},
\end{equation*}
where $j=1,2$ and $\partial\omega_j=\Gamma_{\omega, j}\cup \gamma$.
\end{defi}
Restricting (Eq.~\ref{eq:Psi1_zerleg}), (Eq.~\ref{eq:Psi2_zerleg}) to $I(\gamma)$, we obtain from (Eq.~\ref{eq:neumann}) and $g_{1,\gamma } = g_{2,\gamma } = :g_\gamma $ that
\begin{equation*}
\left( {{\ }^\gamma \Psi_{\omega_1 }^\gamma + {\ }^\gamma \Psi_{\omega_2
}^\gamma } \right)g_\gamma = (- \Psi_{\omega_1 }^f f_1 - \Psi_{\omega_1
}^\Gamma g_{1,\Gamma } - \Psi_{\omega_2 }^f f_2 - \Psi_{\omega_2 }^\Gamma
g_{2,\Gamma })\vert_{I(\gamma)} .
\end{equation*}
Next, we set 
$M: = - ( {{\ }^\gamma \Psi_{\omega_1 }^\gamma + {\ }^\gamma \Psi_{\omega_2
}^\gamma } )$.
and after computing $M^{ - 1}$, we obtain: 
\begin{equation}\label{eq:gamma}
g_\gamma = M^{ - 1}( {\Psi_{\omega_1 }^f f_1 + \Psi_{\omega_1 }^\Gamma
g_{1,\Gamma } + \Psi_{\omega_2 }^f f_2 + \Psi_{\omega_2 }^\Gamma g_{2,\Gamma
} })\vert_{I(\gamma)}.
\end{equation}
\begin{rem}
The inverse matrix $M^{ - 1}$ exists since it is the sum of positive definite matrices
corresponding to the mappings 
${\ }^\gamma \Psi_{\omega_1 }^\gamma, {\ }^\gamma \Psi_{\omega_2}^\gamma$.
\end{rem}
\begin{rem}
Since $g_{\gamma ,i} = u_h (\x_i )$, $i\in I(\gamma)$, we have determined the map $\Phi_\omega $
(it acts on the data $d_{\omega}$ composed by $f_1$, $f_2$, $g_{1,\Gamma
},g_{2,\Gamma }$).
\end{rem}
\begin{rem}
\label{rem:psi}
We have the formula $\Psi_{\omega}(d_\omega)=\Psi_{\omega_1}(d_1)+\Psi_{\omega_2}(d_2)$, where
\begin{equation}
\label{eq:Psi_sum}
\begin{array}{c}
d_\omega = (f_\omega,g_\omega), \quad
d_1 = ( f_1 ,g_{1,\Gamma} ,g_{1,\gamma}), \quad d_2 = ( f_2 ,g_{2,\Gamma
}, g_{2,\gamma } ) , \\ 
g_{1,\gamma } = g_{2,\gamma} = M^{ - 1}({\Psi_{\omega_1 }^f f_1 +
\Psi_{\omega_1 }^\Gamma g_{1,\Gamma } + \Psi_{\omega_2 }^f f_2 +
\Psi_{\omega_2 }^\Gamma g_{2,\Gamma }} )\vert_{I(\gamma)}.
\end{array}   
\end{equation}
Here $(f_{\omega},g_{\omega})$ is build as in (Eq.~\ref{eq:f})-(Eq.~\ref{eq:g}) and 
(Eq.~\ref{eq:dirich}),(Eq.~\ref{eq:contin}) are satisfied.
\end{rem}
\textbf{Conclusion:}\\
Thus, using the given mappings $\Psi_{\omega_1}$, $\Psi_{\omega_2}$, defined on the sons 
$\omega_1, \omega_2 \in T_{\Th}$, we can compute $\Phi_\omega$ and $\Psi_\omega$ for the
father $\omega\in T_{\Th}$. 
\subsection{Building of Matrices $\Psi_{\omega}$ and $\Phi_{\omega}$ from $\Psi_{\omega_1}$ and $\Psi_{\omega_2}$}
\label{sec:build_father}
Let $\omega$, $\omega_1$ where $\omega_2 \in T_{\Th}$ and $\omega_1$, $\omega_2$ are sons of $\omega$. Recall that $\partial\omega_i=\Gamma_{\omega,i}\cup \gamma$. Suppose we have two linear systems of equations for $\omega_1$ and $\omega_2$ which can be written in the block-matrix form:
\begin{equation}
\label{eq:block_equation3}
\left(\begin{array}{cc}
A_{11}^{(i)} & A_{12}^{(i)} \\ 
A_{21}^{(i)} & A_{22}^{(i)}
\end{array}\right)
\left(\begin{array}{c}
\u_1^{(i)} \\ 
\u_2^{(i)}
\end{array}\right)
=
\left(\begin{array}{cc}
F_{11}^{(i)} & F_{12}^{(i)} \\ 
F_{21}^{(i)} & F_{22}^{(i)}
\end{array}\right)
\left(\begin{array}{c}
\c_1^{(i)} \\ 
\c_2^{(i)}
\end{array}\right), \quad i=1,2,
\end{equation}
where $\gamma:=\gamma_{\omega}$,
\begin{equation*}
A_{11}^{(i)}:\RR^{I(\Gamma_{\omega,i})}\rightarrow \RR^{I(\Gamma_{\omega,i})},\quad
A_{12}^{(i)}:\RR^{I(\gamma)}\rightarrow \RR^{I(\Gamma_{\omega,i})},
\end{equation*}
\begin{equation*}
A_{21}^{(i)}:\RR^{I(\Gamma_{\omega,i})}\rightarrow \RR^{I(\gamma)},\quad
A_{22}^{(i)}:\RR^{I(\gamma)}\rightarrow \RR^{I(\gamma)},
\end{equation*}
\begin{equation*}
F_{11}^{(i)}:\RR^{I(\omega_i \setminus \gamma)}\rightarrow \RR^{I(\partial \omega_i)},\quad
F_{12}^{(i)}:\RR^{I(\gamma)}\rightarrow \RR^{I(\partial\omega_i)},
\end{equation*}
\begin{equation*}
F_{21}^{(i)}:\RR^{I(\omega_i \setminus \gamma)}\rightarrow \RR^{I(\gamma)},\quad
F_{22}^{(i)}:\RR^{I(\gamma)}\rightarrow \RR^{I(\gamma)}.
\end{equation*}
Both the equations in (Eq.~\ref{eq:block_equation3}) are analogous to (Eq.~\ref{eq:Psi1_zerleg}) and (Eq.~\ref{eq:Psi2_zerleg}).
Note that $\c_2^{(1)}=\c_2^{(2)}$ and $\u_2^{(1)}=\u_2^{(2)}$ because of the consistency conditions (see (Eq.~\ref{eq:dirich}),(Eq.~\ref{eq:contin})) on the interface $\gamma$. The system of linear equations for $\omega$ be
\begin{equation}
\label{eq:block_equation4}
\left(\begin{array}{ccc}
A_{11}^{(1)} & 0 &A_{12}^{(1)} \\ 
0 & A_{11}^{(2)} &A_{12}^{(2)} \\ 
A_{21}^{(1)} & A_{21}^{(2)} & A_{22}^{(1)}+A_{22}^{(2)} 
\end{array}\right)
\left(\begin{array}{c}
\u_1^{(1)} \\ 
\u_1^{(2)} \\
\u_2^{(1)}  
\end{array}\right)
=
\left(\begin{array}{ccc}
F_{11}^{(1)} & 0 & F_{12}^{(1)} \\ 
0 & F_{11}^{(2)} & F_{12}^{(2)} \\
F_{21}^{(1)} & F_{21}^{(2)} & F_{22}^{(1)}+F_{22}^{(2)}
\end{array}\right)
\left(\begin{array}{c}
\c_1^{(1)} \\ 
\c_1^{(2)} \\
\c_2^{(1)} \\
\end{array}\right).
\end{equation}
See the left matrix in Fig.~\ref{fig:h2hg} in the Appendix.
Using the notation 
\begin{equation*}
\tA_{11}:= 
\left(\begin{array}{cc}
A_{11}^{(1)} & 0 \\ 
0 & A_{11}^{(2)}
\end{array}\right)
,\quad 
\tA_{12}:= 
\left(\begin{array}{c}
A_{12}^{(1)} \\ 
A_{12}^{(2)}
\end{array}\right), 
\end{equation*}

\begin{equation*}
\tA_{21}:=(A_{21}^{(1)}, A_{21}^{(2)}),\quad
\tA_{22}:=A_{22}^{(1)}+A_{22}^{(2)},
\end{equation*}
\begin{equation*}
\tilde{\u}_1:=\left(\begin{array}{c}
\u^{(1)}_1 \\ 
\u^{(2)}_1
\end{array}\right), \quad \tilde{\u}_2:=\u_2^{(1)}=\u_2^{(2)},
\end{equation*}
\begin{equation*}
\tilde{F}_1:=
\left(\begin{array}{ccc}
F_{11}^{(1)} & 0 & F_{12}^{(1)} \\ 
0 & F_{11}^{(2)} & F_{12}^{(2)} 
\end{array}\right), \quad
\tilde{F}_2:=\left(\begin{array}{ccc}
F_{21}^{(1)}, & F_{21}^{(2)}, & F_{22}^{(1)}+F_{22}^{(2)}
\end{array}\right),
\end{equation*}
\begin{equation*}
\tilde{\c}_1:=\left(\begin{array}{c}
\c^{(1)}_1 \\ 
\c^{(2)}_1
\end{array}\right), \quad \tilde{\c}_2:=\c_2^{(1)}=\c_2^{(2)}, \quad
\tilde{\c}:=\left(\begin{array}{c}
\tilde{\c}_1 \\
\tilde{\c}_2
\end{array}\right),
\end{equation*}
the system (Eq.~\ref{eq:block_equation4}) can be rewritten as
\begin{equation}
\label{eq:block_equation5}
\left(\begin{array}{cc}
\tA_{11}^{(i)} & \tA_{12}^{(i)} \\ 
\tA_{21}^{(i)} & \tA_{22}^{(i)}
\end{array}\right)
\left(\begin{array}{c}
\tilde{\u}_1 \\ 
\tilde{\u}_2
\end{array}\right)
=
\left(\begin{array}{c}
\tilde{F}_{1} \\ 
\tilde{F}_{2}
\end{array}\right)
\tilde{\c}.
\end{equation}
The system (Eq.~\ref{eq:block_equation5}), indeed, coincides with (Eq.~\ref{eq:block_equation1}).
After elimination of variables $\u^{(2)}_1$ (on the interface), we obtain the matrices as it shown in  Figures \ref{fig:Hg} and \ref{fig:Hrhs}.
\subsection{Algorithms ``Leaves to Root'' and ``Root to Leaves''}
\label{sec:root2leaves}
The scheme of the recursive process of computing $\Psi_\omega$ and $\Phi_\omega$ from $\Psi_{\omega_1}$ and $\Psi_{\omega_2}$ for all $\omega \in T_{\Th}$ is shown in Fig. \ref{fig:HDD_th2} (left). We call this process ``\textbf{Leaves to Root}'':
\begin{enumerate}
\item Compute $\Psi^f_\omega \in \RR^{3\times 3}$ and $\Psi^g_\omega \in \RR^{3\times 3}$ on all leaves of $T_{\Th}$ (triangles of $\mathcal{T}_h$) by (Eq.~\ref{eq:coefA}) and (Eq.~\ref{eq:coefF}).
\item Compute recursive from leaves to root $\Phi_\omega$ and $\Psi_\omega $ from $%
\Psi_{\omega_1 } ,\Psi_{\omega_2 } $. Store $\Phi_{\omega}$ and delete $\Psi_{\omega_1 } ,\Psi_{\omega_2 }$.
\item Stop if $\omega=\Omega$.
\end{enumerate}
\begin{rem}
The result of this algorithm will be a collection of mappings $\{\Phi_\omega
:\omega \in T_{\mathcal{T}_h}\}$. The mappings $\Psi_\omega $, $\omega\in T_{\Th}$, are only of auxiliary purpose and need not stored. 
\end{rem}
The algorithm which applies the mappings $\Phi_{\omega}=(\Phi^g_{\omega}, \Phi^f_{\omega})$ to compute the solution we call ``\textbf{Root to Leaves}''. This algorithm starts from the root and ends on the leaves. Figure \ref{fig:HDD_th2} (right) presents the scheme of this algorithm. \\
Let the data $d_\omega = (f_\omega ,g_\omega )$, $\omega =\Omega$, be given. We can then compute the solution $u_h$ of the initial problem as follows.\\
The Algorithm ``\textbf{Root to Leaves}'':
\begin{enumerate}
\item Start with $\omega=\Omega$.
\item Given $d_\omega = (f_\omega ,g_\omega )$, compute the solution $u_h $
on the interior boundary $\gamma _\omega $ by $\Phi_\omega (d_\omega )$.
\item Build the data $d_{\omega_1}=(f_{\omega_1 } ,g_{\omega_1 } )$, $%
d_{\omega_2}=(f_{\omega_2}, g_{\omega_2})$ from $d_\omega = (f_\omega
,g_\omega )$ and \linebreak $g_{\gamma_{\omega}} : = \Phi_\omega (d_\omega )$.
\item Repeat the same for the sons of $\omega_1$ and $\omega_2$.
\item End if $\omega$ does not contain internal nodes.
\end{enumerate}
Since $u_h(\x_i)=g_{\gamma,i}$, the set of values $(g_{\gamma_{\omega}})$, for all $\omega \in T_{\Th}$, results the
solution of the initial problem (Eq.~\ref{eq:elliptic}) in the whole domain $\Omega$. 

\subsection{Multiple scales}
\label{ssec:multiscales}
Let $h$ and $H$ be fine and coarse meshes, used for discretization of Eq.~\ref{eq:elliptic_det}.
The subscript ${}_h$ near the operator or function means that this operator or function was discretized on a mesh with the step size $h$. 
Let $n_h$ and $n_H$ be the numbers of degrees of freedom on a fine grid and on a coarse grid. For instance, if the right-hand side is smooth, we may use a coarser mesh for it (e.g., operator $\F_H$ and $\G_h$). So, the matrices $\Psi^f and \Phi^f$ will be much smaller. For discretising the diffusion coefficient and the Dirichlet b.c. we use a fine scale $h$ (see more in \cite{MYPHD}).
\begin{lemma}
The complexities of the one-grid version and two-grid version of HDD are 
\begin{equation*}
\O(n_h\log^3 n_h)\quad\text{and} \quad
\O(\sqrt{n_hn_H}\log^3 \sqrt{n_hn_H}),
\quad \text{respectively}.
\end{equation*}
The storage requirements of the one-grid version and two-grid version of HDD are \begin{equation*}
\O(n_h\log^2n_h) \quad\text{and}\quad \O(\sqrt{n_hn_H}\log^2 \sqrt{n_hn_H}),
\quad \text{respectively}.
\end{equation*}
\end{lemma}
\Proof see \cite{MYPHD, FlorianDiss} or Ch. 12 in \cite{HackHMBook}.

\section{Hierarchical matrix approximation}
\label{sec:Hmatrices}

\begin{figure}[htbp]
\centerline{\includegraphics[width=0.3\textwidth]{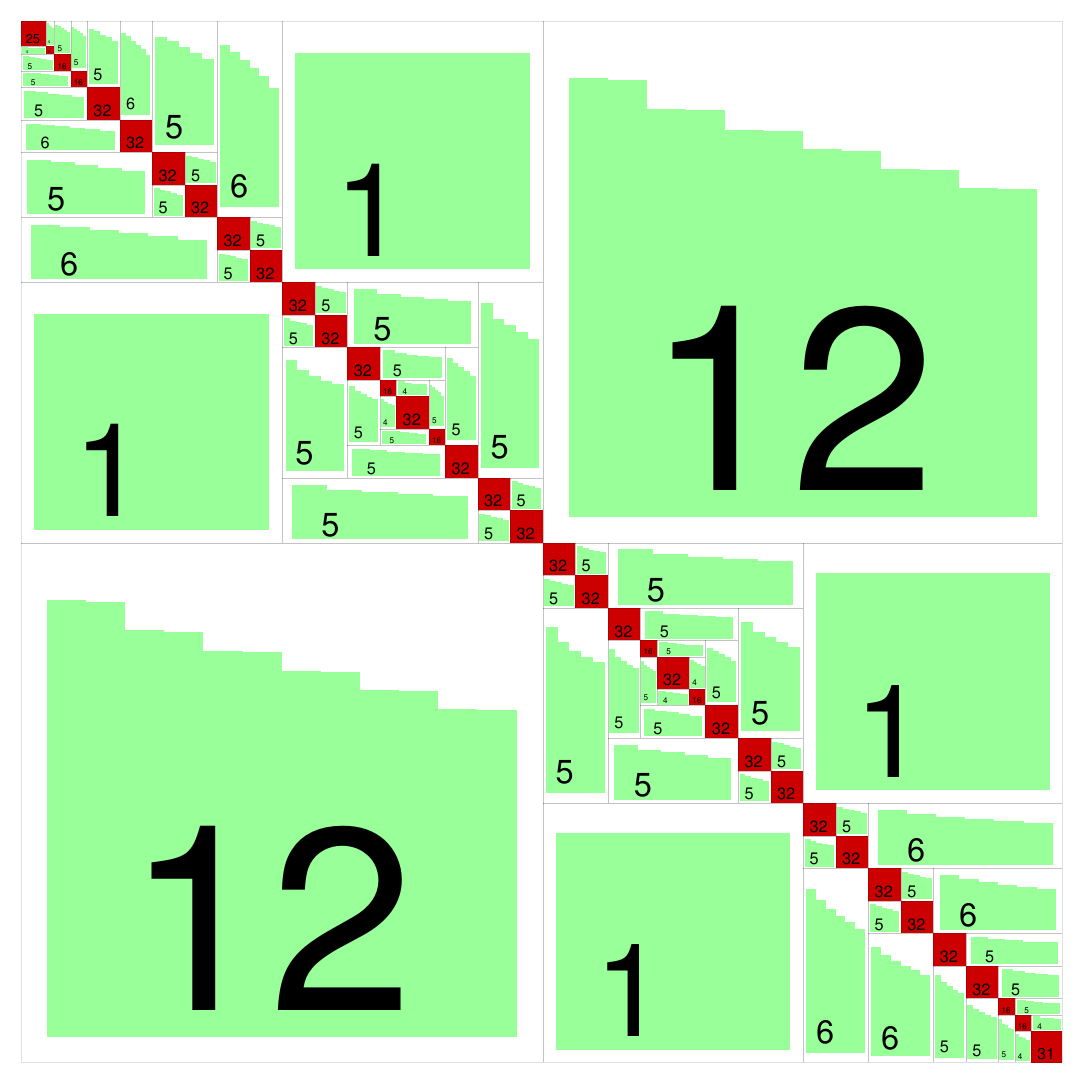}}
\caption{An $\H$-matrix approximation to $(\Psi_{\omega}^g)^{\H} \in \RR^{I\times I}$, $I:=I(\partial \omega )$. The dark (red) blocks are dense matrices and  grey (green) blocks are low-rank matrices. The numbers inside the blocks indicate the ranks of these blocks.}
\label{fig:Hg}
\end{figure}

\begin{figure}[htbp]
\centerline{\includegraphics[width=0.99\textwidth]{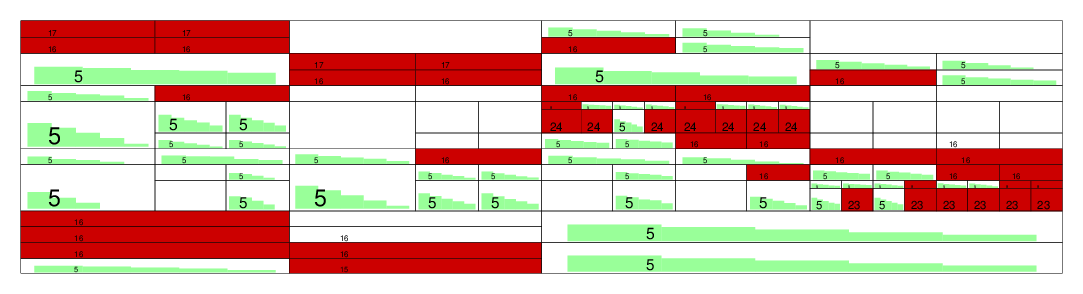}}
\caption{An $\H$-matrix approximation to $(\Psi_{\omega}^f)^{\H} \in \RR^{I\times J}$, $I:=I(\partial \omega)$, $J:=J(\omega)$, $\vert I\vert=256$, $\vert I\vert=4225$. The dark (or red) blocks indicate the dense matrices and the grey (green) blocks indicate the
rank-$k$ matrices; the number inside each block is its rank. The steps inside the blocks show the decay of
the singular values in log scale. The white blocks are empty. }
\label{fig:Hrhs}
\end{figure}

The mappings $\Psi_{\omega}$
and $\Phi_{\omega}$ correspond to dense matrices, and, therefore, require quadratic storage and quadratic or cubic arithmetic cost.
Both these mappings (matrices) $\Psi_{\omega}$ (see Fig.~\ref{fig:Hg}) and
$\Phi_{\omega}$ (see Fig.~\ref{fig:Hrhs}) can be approximated in the $\mathcal{H}$-matrix format. Additionally, all necessary computational steps can be performed within the hierarchical matrix format with a log-linear cost.

The matrices 
$\Phi^g_{\omega}:\RR^{I(\partial\omega)}\rightarrow \RR^{I(\gamma_{\omega})}$, 
$\Phi^f_{\omega}:\RR^{I(\omega)}\rightarrow \RR^{I(\gamma_{\omega})}$,
$\Psi^f_{\omega}:\RR^{ I(\omega) }\rightarrow  \RR^{ I(\partial\omega) }$,
are rectangular. 
The matrix $\Psi^g_{\omega}:\RR^{ I(\partial\omega) }\rightarrow  \RR^{ I(\partial\omega) }$ is quadratic. 
%
%
%
%

The hierarchical matrices ($\H$-matrices) have been used in a wide range of applications since their introduction in 1999 by Hackbusch \cite{Part1}. They provide a format for the data-sparse representation of fully-populated matrices. 
The complexity of the $\mathcal{H}$-matrix addition, multiplication, Schur complement and inversion is $\mathcal{O}(k^2n\log^q n)$, $q=1,2$. See more details about $\H$-matrices in \cite{HackHMBook,Part1,part2,GH03,PhD,Litv_Genton2019,Litv_HLIBCov2020}.
%
%
In \cite{BoxLU} authors prove the existence of an $\H$-matrix approximation of the inverse (Assumption~2) and of the Schur complement (Theorem~1).

The following proposition follows from Theorem~1 (\cite{BoxLU}) and \cite{weak}. 
\begin{prop}
The matrices $\Psi^g_{\omega}\in \RR^{I(\partial\omega)\times I(\partial\omega)}$, $\Psi^f_{\omega} \in \RR^{I(\partial\omega)\times I(\omega)}$, $\Phi^f_{\omega}\in \RR^{I(\partial\omega)\times I(\omega)}$ and $\Phi^g_{\omega}\in \RR^{I(\gamma_\omega)\times I(\partial\omega)}$ for all $\omega \in T_{\Th}$ can be effectively approximated by $\H$-matrices.
\end{prop}
For more details and complexity estimates see \cite{MYPHD}.


\section{Fast Evaluation of Functionals} 
\label{sec:functional}
\begin{figure}[ht!]
\centerline{\includegraphics[width=0.25\textwidth]{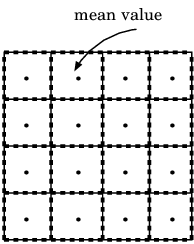}\quad \quad
\includegraphics[width=0.45\textwidth]{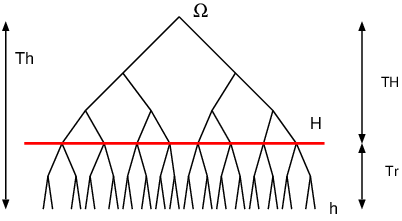}}
\caption{(left) HDD computes the solution on a coarse mesh (or the interface $I(\Sigma)$) and the mean value inside of each cell. (right) Algorithm ``Leaves to Root" goes through the whole tree, but ``Root to leaves" starts in the root, goes through a subtree $T_{\TH}$ and terminates on a coarse level with the mesh size $H$ (marked with the red horizontal line). After that the mean value inside of each coarse cell (of size $H\times H$) is computed.}
\label{fig:sol_cell}
\end{figure}

%
In this section we describe how to use $\Phi^f_{\omega}$ and $\Phi^g_{\omega}$ for building different linear functionals of the solution. Indeed, the functional $\lambda$ is determined in the same way as $\Psi_{\omega}$.\\
Below we list some examples of problem settings which can be solved by using linear functionals.
\begin{ex}If the solution $u$ in a subdomain $\omega \in T_{\Th}$ is known, the mean value $\mu(\omega)$ can be computed by the following formula
\begin{equation}
\label{eq:mean}
\overline{u}_{\omega}:=\mu(\omega)=\displaystyle{\frac{\int_{\omega}{u(\x)d\x}}{\vert \omega \vert}=\frac{\sum_{t \in \Th(\omega)}{\frac{\vert t \vert}{3}(u_1+u_2+u_3)}}{\vert \omega \vert}},
\end{equation}
where $u$ is affine on each triangle $t$ with values $u_1$, $u_2$, $u_3$ at the three corners and $\Th(\omega)$ is the collection of all triangles in $\omega$. If the solution $u$ is unknown, we would like to have a linear functional $\l_{\omega}(f,g), \omega \in T_{\Th}$, which computes the mean value $\mu_{\omega}$ of the solution in $\omega$. \\
\end{ex}
\begin{ex}
\label{ex:MeanValue}
Let us assume that $\Omega$ is decomposed into $p=16$ subdomains $\Omega=\bigcup_{i=1}^{p}\Omega_i$. Sometimes these sub-domains $\Omega_i$ are called \texttt{cells}.
The set of nodal points on the interface is denoted by $I_{\Sigma}$. HDD can compute the solution on the interface $I_{\Sigma}$ and the mean value over each $\Omega_i$, $i=1,...,p$, (see Fig.~\ref{fig:sol_cell}).
\end{ex}
\begin{ex}
To compute the FE solution $u_h(\x_i)$ in a fixed nodal point $\x_i \in \Omega$, i.e., to define how the solution $u_h(\x_i)$ depends on the given FE Dirichlet data $g_h\in \RR^{I(\partial\Omega)}$ and the FE right-hand side $f_h\in \RR^{I(\Omega)}$.
\end{ex}
%
\subsection{Computing the mean value in all subdomains $\omega \in T_{\Th}$}
Let $\omega$ be the father node, $\omega_1$ the left son, and $\omega_2$ the right son.
Then $\omega=\omega_1 \cup \omega_2$, $\omega_1\cap \omega_2\neq \varnothing$, with $\omega$, $\omega_1$, $\omega_2 \in T_{\Th}$.
To simplify the notation, we will write $d_i:=d_{\omega_i}$ and $(f_i,g_i)$ instead of $(f_{\omega_i},g_{\omega_i})$, $i=1,2$ (see also Fig.~\ref{fig:ex_points}).
Recall the following notation (see (Eq.~\ref{eq:d_components}), (Eq.~\ref{eq:g_components})):
\begin{equation*}
\Gamma:=\partial\omega, \quad \Gamma_{\omega,1}:=\partial\omega \cap \omega_1,\quad  \Gamma_{\omega,2}:=\partial\omega \cap \omega_2, \quad \text{then}
\end{equation*}
\begin{equation}
\label{eq:d1_d2}
d_1 = (f_1, g_1) = (f_1, g_{1,\Gamma }, g_{1,\gamma } ),\quad d_2 = (f_2, g_2) = (f_2, g_{2,\Gamma } ,g_{2,\gamma}), \text{ where }
\end{equation}
\begin{equation}
\label{eq:g1_g2}
\begin{array}{c}
g_{1,\Gamma } := (g_1 )\vert_{\Gamma _{\omega,1}}, \quad g_{1,\gamma } := (g_1)\vert_{\gamma},\\
g_{2,\Gamma } := (g_2 )\vert_{\Gamma _{\omega,2}}, \quad g_{2,\gamma } := (g_2)\vert_{\gamma}.
\end{array}
\end{equation}
We consider a linear functional $\l_{\omega}$ with the following properties:
\begin{equation}
\label{eq:prop1}
\l_{\omega}(d_{\omega})=(\l^g_{\omega},g_{\omega})+(\l^f_{\omega},f_{\omega}),
\end{equation}
\begin{equation}
\label{eq:prop2}
\l_{\omega}(d_{\omega})=c_1\l_{\omega_1}(d_{\omega_1})+c_2\l_{\omega_2}(d_{\omega_2}),
\end{equation}
where $\l^g_{\omega}:\RR^{I(\partial\omega)}\rightarrow \RR$, $\l^f_{\omega}:\RR^{I(\omega)}\rightarrow \RR$, $c_1$, $c_2$ two constants, and $(\cdot,\cdot)$ the scalar product of two vectors.
%
%
%
\begin{defi}
\label{eq:vextension}
Let $\omega_1\subset \omega$, $\l^f_{\omega_1}: \RR^{I(\omega_1)} \rightarrow \RR$. a) We define the following extension of $\l^f_{\omega_1}$ 
\begin{equation*}
(\l_{\omega_1}^f\vert^{\omega})_i:=
\left\{ {{\begin{array}{*{20}ll} 
(\l_{\omega_1}^f)_i & \text{for } i \in I(\omega_1), \hfill \\ 
0 & \text{for } i \in I(\omega \setminus \omega_1 ), \hfill \\ \end{array} }} \right.
\end{equation*}
where $(\l_{\omega_1}^f\vert^{\omega}):\RR^{I(\omega)}\rightarrow \RR$.
b) The extension of the functional $\l_{1,\Gamma}^g: \RR^{I(\Gamma_{\omega,1})}\rightarrow \RR$ is defined as 
\begin{equation*}
(\l_{1,\Gamma}^g\vert^{\Gamma})_i:=
\left\{ {{\begin{array}{*{20}ll} 
(\l_{1,\Gamma}^g)_i & \text{for } i \in I(\Gamma_{\omega,1}), \hfill \\ 
0 & \text{for } i \in I(\Gamma \setminus \Gamma_{\omega,1}), \hfill \\ \end{array} }} \right.
\end{equation*}
where $(\l_{1,\Gamma}^g\vert^{\Gamma}):\RR^{I(\Gamma)}\rightarrow \RR$.
\end{defi}
%
%
\begin{defi}
Using (Eq.~\ref{eq:g1_g2}), we obtain the following decompositions\\ $\l_{\omega_1}^g=(\l_{1,\Gamma}^g,\l_{1,\gamma}^g)$ and $\l_{\omega_2}^g=(\l_{2,\Gamma}^g,\l_{2,\gamma}^g)$, where 
$\l_{1,\Gamma}^g:\RR^{I(\Gamma_{\omega,1})}\rightarrow \RR$,
$\l_{1,\gamma}^g:\RR^{I(\gamma)}\rightarrow \RR$,\\
$\l_{2,\Gamma}^g:\RR^{I(\Gamma_{\omega,2})}\rightarrow \RR$,
$\l_{2,\gamma}^g:\RR^{I(\gamma)}\rightarrow \RR$.
\end{defi}
\begin{lemma}
\label{lem:Func1}
Let $\l_{\omega}(d_{\omega})$ satisfy (Eq.~\ref{eq:prop1}) and (Eq.~\ref{eq:prop2}) with $\omega=\omega_1 \cup \omega_2$. Let $\l^g_{\omega_1}$, $\l^g_{\omega_2}$, $\l^f_{\omega_1}$ and $\l^f_{\omega_2}$ be the vectors for the representation of the functionals $\l_{\omega_1}(d_{\omega_1})$ and $\l_{\omega_2}(d_{\omega_2})$. Then the vectors $\l_{\omega}^f$, $\l_{\omega}^g$ for the representation
\begin{equation}
\label{eq:fformula}
\l_{\omega}(d_{\omega})=(\l_{\omega}^f,f_{\omega})+(\l_{\omega}^g,g_{\omega}), 
\text{ where } f_{\omega}\in \RR^{I(\omega) }, g_{\omega} \in \RR^{I(\partial\omega)},\\
\end{equation}
are given by
\begin{equation*}
\l_{\omega}^f=\tilde{\l}^{f}_{\omega}+(\Phi^f_{\omega})^T \l^g_{\gamma},
\end{equation*}
\begin{equation*}
\l_{\omega}^g=\tilde{\l}^{g}_{\omega}+(\Phi^g_{\omega})^T \l^g_{\gamma},
\end{equation*}
\begin{equation}
\label{eq:funct_f1}
\tilde{\l}^f_{\omega}:=c_1\l_{\omega_1}^f\vert^{\omega}+c_2\l_{\omega_2}^f\vert^{\omega}, 
\end{equation}
\begin{equation}
\label{eq:funct_g1}
\tilde{\l}^g_{\Gamma}:=c_1\l^g_{1,\Gamma}\vert^{\Gamma}+c_2\l^g_{2,\Gamma}\vert^{\Gamma},
\end{equation}
\begin{equation}
\label{eq:funct_g2}
\l^g_{\gamma}=c_1\l^g_{1,\gamma}+c_2\l^g_{2,\gamma}.
\end{equation}
\end{lemma}
\Proof Let $d_{\omega_1}$, $d_{\omega_2}$ be the given data and
$\l_{\omega_1}$ and $\l_{\omega_2}$ be the given functionals. Then the functional $\l_{\omega}$ satisfies
\begin{align*}
\l_{\omega}(d_{\omega})&\overset{(Eq.~\ref{eq:prop2})}{=}c_1\l_{\omega_1}(d_{\omega_1})+c_2\l_{\omega_2}(d_{\omega_2})\\
&\overset{(Eq.~\ref{eq:prop1})}{=}c_1((\l_{\omega_1}^f,f_1)+(\l_{\omega_1}^g,g_1))+c_2((\l_{\omega_2}^f,f_2)+(\l_{\omega_2}^g,g_2)).
\end{align*}
Using the decomposition (Eq.~\ref{eq:d1_d2}), we obtain
\begin{align}
\label{eq:lambda1}
\l_{\omega}(d_{\omega})&=c_1(\l_{\omega_1}^f,f_1)+c_2(\l_{\omega_2}^f,f_2)+c_1((\l^g_{1,\Gamma},g_{1,\Gamma})+(\l^g_{1,\gamma},g_{1,\gamma}))\\
&+c_2((\l^g_{2,\Gamma},g_{2,\Gamma})+(\l^g_{2,\gamma},g_{2,\gamma})).
\end{align}
The consistency of the solution implies $g_{1,\gamma}=g_{2,\gamma}=:g_{\gamma}$. From the Definition \ref{eq:vextension} follows
\begin{equation*}
(\l_{\omega_1}^f,f_1)=(\l_{\omega_1}^f\vert^{\omega},f_{\omega}),\quad (\l_{\omega_2}^f,f_2)=(\l_{\omega_2}^f\vert^{\omega},f_{\omega}),
\end{equation*}
\begin{equation*}
(\l^g_{1,\Gamma},g_{1,\Gamma})=(\l^g_{1,\Gamma}\vert^{\Gamma},g_{\omega}), \quad
(\l^g_{2,\Gamma},g_{2,\Gamma})=(\l^g_{2,\Gamma}\vert^{\Gamma},g_{\omega}).
\end{equation*}
Then, we substitute the last expressions in (Eq.~\ref{eq:lambda1}) to obtain
\begin{equation}
\label{eq:Lomega}
\l_{\omega}(d_{\omega})=(c_1\l_{\omega_1}^f\vert^{\omega}+c_2\l_{\omega_2}^f\vert^{\omega},f_{\omega})+(c_1\l^g_{1,\Gamma}\vert^{\Gamma}+c_2\l^g_{2,\Gamma}\vert^{\Gamma},g_{\omega})
\end{equation}
\begin{equation*}
+(c_1\l^g_{1,\gamma}+c_2\l^g_{2,\gamma},g_{\gamma}).
\end{equation*}
%
Set $\tilde{\l}^f_{\omega}:=c_1\l_{\omega_1}^f\vert^{\omega}+c_2\l_{\omega_2}^f\vert^{\omega}$, 
$\tilde{\l}^g_{\Gamma}:=c_1\l^g_{1,\Gamma}\vert^{\Gamma}+c_2\l^g_{2,\Gamma}\vert^{\Gamma}$and $\l_{\gamma}^g:=c_1\l_{1,\gamma}^g+c_2\l_{2,\gamma}^g$.
\\
From the algorithm ``Root to Leaves'' we know that 
\begin{equation}
\label{eq:g_gamma}
g_{\gamma}=\Phi_{\omega}(d_{\omega})=\Phi^g_{\omega}\cdot g_{\omega}+\Phi^f_{\omega}\cdot f_{\omega}. 
\end{equation}

Substituting $g_{\gamma}$ from (Eq.~\ref{eq:g_gamma}) in (Eq.~\ref{eq:Lomega}), we obtain 
\begin{align*}
\l_{\omega}(d_{\omega})&=(\tilde{\l}^{f}_{\omega},f_{\omega})+(\tilde{\l}^g_{\omega},g_{\omega})+(\l^g_{\gamma},\Phi^g_{\omega} g_{\omega}+\Phi^f_{\omega} f_{\omega})\\
&=(\tilde{\l}^{f}_{\omega}+(\Phi^f_{\omega})^T \l^g_{\gamma},f_{\omega})+(\tilde{\l}^{g}_{\omega}+(\Phi^g_{\omega})^T \l^g_{\gamma},g_{\omega}).
\end{align*}
We define $\l^f_{\omega}:=\tilde{\l}^{f}_{\omega}+(\Phi^f_{\omega})^T \l^g_{\gamma}$ and $\l^g_{\omega}:=\tilde{\l}^{g}_{\omega}+(\Phi^g_{\omega})^T \l^g_{\gamma}$ and obtain
\begin{equation}
\l_{\omega}(d_{\omega})=(\l_{\omega}^f,f_{\omega})+(\l_{\omega}^g,g_{\omega}).
\end{equation}
\begin{flushright}
$\blacksquare$
\end{flushright}
\begin{ex}
Lemma \ref{lem:Func1} with $c_1=\frac{\vert \omega_1 \vert}{\vert \omega \vert}$, $c_2=\frac{\vert \omega_2 \vert}{\vert \omega \vert}$ can be used to compute the mean values in all subdomains $\omega \in T_{\Th}$.
\end{ex}
%
%

\subsection{Algorithms for computing the mean values}
\label{sec:MeanValue}

%
%
Below we describe two algorithms which are required for computing the mean value in all subdomains $\omega \in T_{\Th}$. These algorithms compute vectors $\l^g_{\omega}$ and $\l^f_{\omega}$ respectively.\\

The initialisation is $\l^g_{\omega}:=(\frac{1}{3},\frac{1}{3},\frac{1}{3})$, $\l^f_{\omega}:=(0,0,0)$ for all leaves of $T_{\Th}$. Let us denote $\l^g_{1}:=\l^g_{\omega_1}$, $\l^g_{2}:=\l^g_{\omega_2}$. The algorithms for building $\l^g_{\omega}$ and $\l^f_{\omega}$
for all $\omega \in T_{\Th}$, which have internal nodes, are the following:
\begin{algo}
   \label{alg:functional}
    \begin{algorithmic}
      \STATE{(Building of $\l^g_{\omega}$)}
      \STATE{\textbf{build$\_$functional$\_$g}($\l^g_{1}$, $\l^g_{2}$, $\Phi^g_{\omega}$,...) }
      \STATE \textbf{begin}
      \STATE $\quad$ \textbf{allocate memory for} $\l^g_{\omega}$;
      \STATE $\quad$ \textbf{for all} $i \in I(\Gamma_{\omega,1})$ \textbf{do}
        \STATE $\quad\quad$ $\l^g_{\omega}[i]+=c_1\l_1^g[i]$;
      \STATE $\quad$ \textbf{for all} $i \in I(\Gamma_{\omega,2})$ \textbf{do}
        \STATE $\quad\quad$ $\l^g_{\omega}[i]+=c_2\l_2^g[i]$;
      \STATE $\quad$ \textbf{for all} $i \in I(\gamma)$ \textbf{do}
        \STATE $\quad\quad$ $z[i]=c_1\l_1^g[i]+c_2\l_2^g[i]$;
      \STATE $\quad$ $v:=(\Phi^g_{\omega})^T\cdot z$; 
      \STATE $\quad$ \textbf{for all} $i \in I(\partial\omega)$ \textbf{do}
      \STATE $\quad\quad$ $\l^g_{\omega}[i]:=\l^g_{\omega}[i]+v[i]$; 
      \STATE $\quad$ \textbf{return $\l^g_{\omega}$}; 
      \STATE \textbf{end}; 
   \end{algorithmic}
\end{algo}
Let $\l^f_{1}:=\l^f_{\omega_1}$, $\l^f_{2}:=\l^f_{\omega_2}$.
\begin{algo}
   \label{alg:functional_f}
    \begin{algorithmic}
      \STATE{(Building of $\l^f_{\omega}$)}
      \STATE{\textbf{build$\_$functional$\_$f}($\l^f_1$, $\l^f_2$, $\Phi^f_{\omega}$,...) }
      \STATE \textbf{begin}
      \STATE{ $\quad$ \textbf{for all} $i \in I(\omega_1\backslash \gamma)$ \textbf{do}}
        \STATE $\quad\quad$ $\l^f_{\omega}[i]+=c_1\l_1^f[i]$;
      \STATE{ $\quad$ \textbf{for all} $i \in I(\omega_2\backslash \gamma)$ \textbf{do}}
        \STATE $\quad\quad$ $\l^f_{\omega}[i]+=c_2\l_2^f[i]$;
      \STATE{ $\quad$ \textbf{for all} $i \in I(\gamma)$ \textbf{do}}
        \STATE $\quad\quad$ $z[i]=c_1\l_1^f[i]+c_2\l_2^f[i]$;
      \STATE $\quad$ $v:=(\Phi^f_{\omega})^T\cdot z$; 
      \STATE $\quad$ \textbf{for all} $i \in I(\omega)$ \textbf{do}
      \STATE $\quad\quad$ $\l^f_{\omega}[i]:=\l^f_{\omega}[i]+v[i]$; 
      \STATE $\quad$ \textbf{ return $\l^f_{\omega}$}; 
      \STATE \textbf{end}; 
   \end{algorithmic}
\end{algo}
\begin{rem}
a) If only the functionals $\l_{\omega}$, $\omega \in T_{\Th}$, are of interest, the maps $\Phi_{\omega}$ do not need to be stored.\\
b) For functionals with local support in some $\omega_0 \in T_{\Th}$, it suffices that $\Phi_{\omega}$ is given for all $\omega\in T_{\Th}$ with $\omega \supset \omega_0$, while $\l_{\omega_0}(d_{\omega_0})$ is associated with $\omega_0 \in T_{\Th}$. The computation of $\Lambda(u_h)=\lambda(d)$ starts with the recursive evaluation of $\Phi_{\omega}$ for all $\omega\supset \omega_0$. Then the data $d_{\omega_0}$ are available and $\lambda_{\omega_0}$ can be applied.
\end{rem}

\subsection{Solution in a subdomain} 
\label{sec:SolSubdomain}
\begin{figure}[ht!]
\centerline{\includegraphics[width=0.2\textwidth]{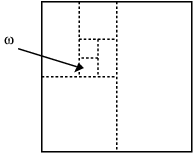}\quad
{\includegraphics[width=0.2\textwidth]{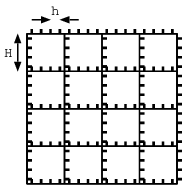}}}
\caption{(left) The solution in a subdomain $\omega \in T_{\Th}$ is required. HDD computes subsequently the solution only on the dotted lines and then only in $\omega$; (right) The coarse $H$ and the fine $h$ scales.}
\label{fig:sol_point}
\end{figure}

Suppose that the solution is only required in a small subdomain $\omega \in T_{\Th}$ (Fig. \ref{fig:sol_point}, left). For this purpose, the HDD method requires less computational resources as usual. The algorithm ``Leaves to Root'' is performed completely, but the algorithm ``Root to Leaves'' computes the solution only on the internal boundaries (dotted lines) which are necessary for computing the solution in $\omega$. The storage requirements are also significantly reduced. We only store the mappings $\Phi^f_{\omega}$ and $\Phi^g_{\omega}$ for all $\omega \in T_{\Th}$ that belong to the path from the root of $T_{\Th}$ to $\omega$. The storage requirement is $\O(n_h\log n_h)$, where $n_h$ is the number of degrees of freedom in $\Omega$. The computational cost of the ``Root to Leaves'' is $\O(n_h\log^2 n_h)$. If the right-hand side is smooth, it can be discretized (defined) only on a coarse mesh $\TH$ (see Fig.~\ref{fig:sol_point}, right). About the interpolation and restriction operators read in \cite{MYPHD,Hackbusch2012,FlorianDiss,HackHMBook}. 
%
%

\section{Numerics}
\label{sec:numeric}
In \cite{MYPHD, FlorianDiss, Hackbusch2012}, the HDD method was compared with the preconditioned conjugate gradient (PCG) method, with the hierarchical ($\mathcal{H}$)- Cholesky method and the direct full $\H$-matrix inverse. 
A cheap $\H$-matrix approximation of the inverse, computed from the $\H$-Cholesky factors, was used as a preconditioner.

Some experiments were performed with two meshes -  a coarse for the right-hand side and a fine for the diffusion coefficient.
Technical details and implementation of the HDD method can be found in \cite{HDDdoc, FlorianDiss, MYPHD}.
The data misfit and the likelihood function in the Bayesian-like approach were computed in 
\cite{rosic2012sampling,rosic2013parameter,hermann2016inverse,matthies2016parameter,pajonk2012deterministic,rosic2011direct}.
In the following numerical experiments we compare the computational time and memory requirement of HDD with the times and memory requirements of the $\H$-matrix inverse and the $\H$-Cholesky factorisation.

\begin{figure}[htbp]
\centerline{
\includegraphics[width=0.49\textwidth]{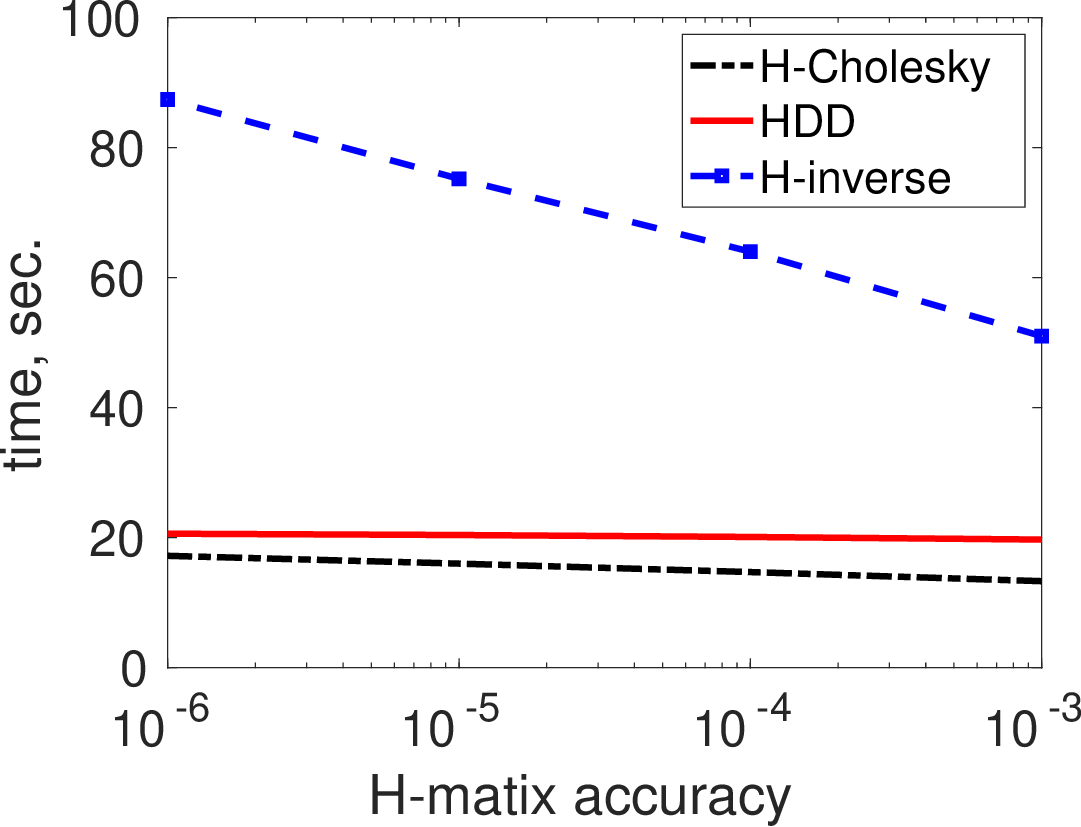}
\includegraphics[width=0.49\textwidth]{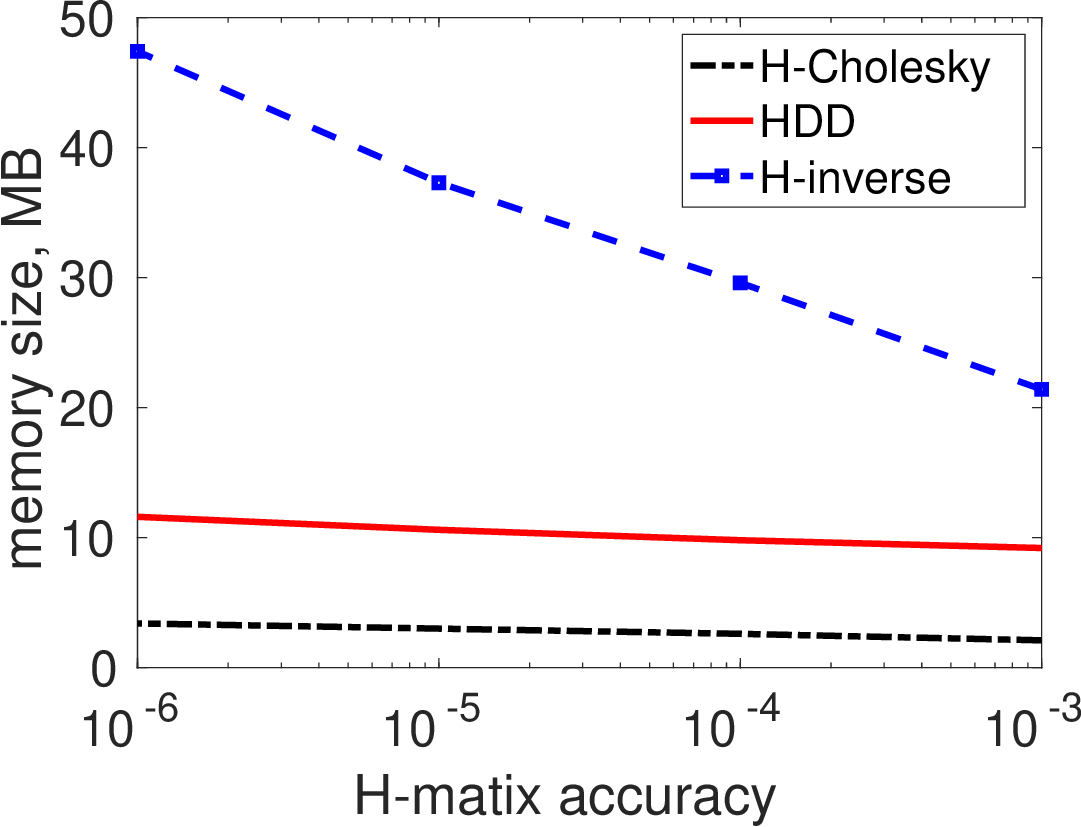}
}
\caption{Comparison of the HDD, $\H$-Cholesky factorisation, and $\H$-matrix inverse. (left) Dependence of the computing time (in sec.) and (right)  memory requirements (in MB) on the $\H$-matrix accuracy, $n=129^2$ dofs.}
\label{fig:time_size}
\end{figure}

We consider the problem as in Eq.~\ref{eq:elliptic} with fixed $\bZ$. The computational domain is a unit square $\Omega=[0,1]^2$, the diffusion coefficient is $\kappa(x,y)=1+0.5\cdot \sin(50x)\sin(50y)$. Note, that HDD does not require an axes parallel triangulation. All numerical experiments were performed on a usual notebook.
Figure~\ref{fig:time_size} demonstrates the dependence of computing times (left) and memory requirements (right) on the $\H$-matrix accuracy (see the adaptive rank arithmetic in \cite{HackHMBook}). One can see that the computational time and storage requirement of the $\H$-Cholesky factorisation are the best. The HDD method shows a slightly larger time and storage than the $\H$-Cholesky factorisation (due to some overhead) and is better than the direct $\H$-matrix inverse. Note that HDD computes more details about the operator and the solution than the $\H$-Cholesky factorisation. 

If the right-hand side is smooth,  we can discretise it on a mesh with the mesh size, for instance, $H=2h$. As result, the matrices $\Phi^g$ and $\Psi^g$ will stay the same, but $\Phi^f$ and $\Psi^f$ will be smaller. See more for the restriction and prolongation operators in \cite{MYPHD}.

\begin{figure}[htbp]
\centerline{\includegraphics[width=1.5in,height=1.5in]{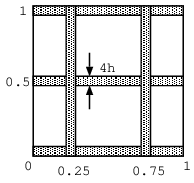}}
\caption{Model domain $\Omega=[0,1]^2$. The diffusion coefficient is very small ($\alpha=10^{-5}$) inside the grey areas and large in white subdomains ($\beta=1$).}
\label{fig:simpleskin}
\end{figure}

In the next example we consider again the problem as in Eq.~\ref{eq:elliptic}. The parameter $\bZ$ is fixed and $\kappa$ is a jumping coefficient as in Fig.~\ref{fig:simpleskin} with $\alpha=10^{-5}$ and $\beta=1$.
Such problems appear in the material sciences and in medicine (the, so-called, skin problem). 

\begin{figure}[htbp]
\centerline{
\includegraphics[width=0.49\textwidth]{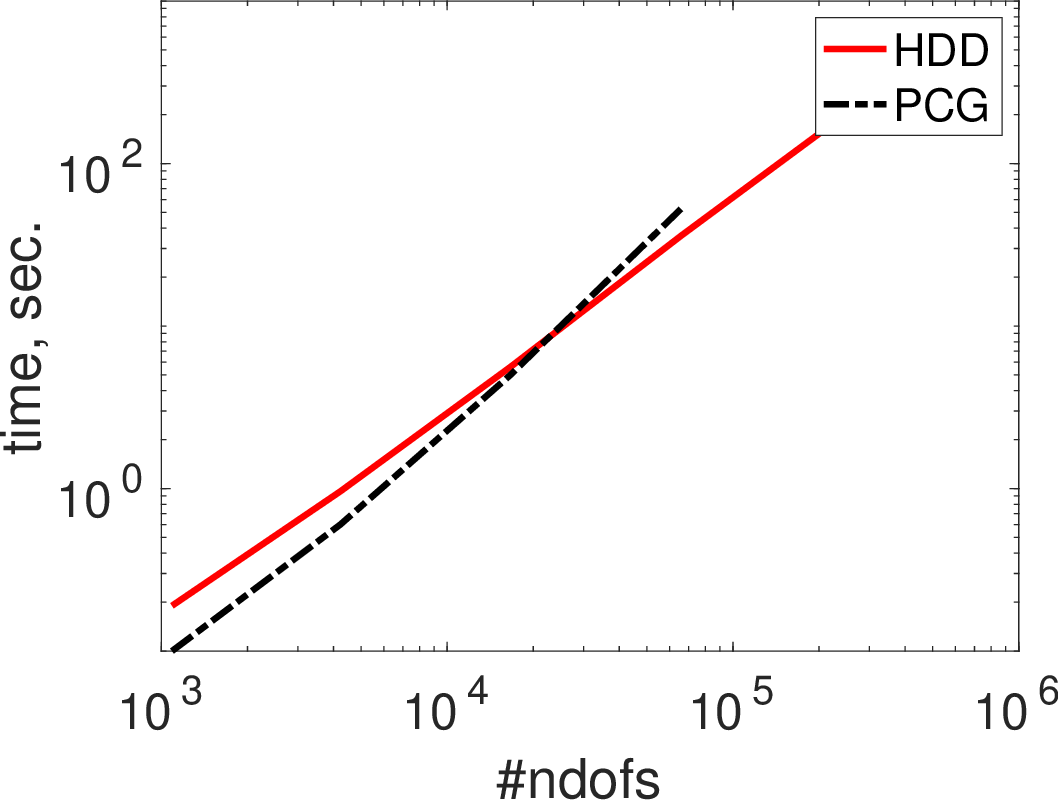}
}
\caption{HDD and PCG computing times vs. $n$.  The accuracy in each $\H$-matrix subblock is $10^{-8}$, the PCG stopping criteria $\varepsilon_{cg}=10^{-8}$, $\frac{H}{h}=2$ as in Sec.~\ref{ssec:multiscales}.}
\label{fig:hdd_pcg_time}
\end{figure}

Figure~\ref{fig:hdd_pcg_time} compares the computational times of HDD and PCG methods. 
The PCG time includes the time needed for: 
(a) computing the stiffness matrix $A$ in the $\H$-matrix format;
(b) computing the $\H$-Cholesky decomposition of $A$ (used as a preconditioner);
(c) PCG iterations.\\ 
For example, for $n=66049$, the PCG time is $53=38.2+11.4+3.4$ (sec.).
Note that for $n\approx 263000$ dofs there is not enough memory to compute the stiffness matrix $A$ and perform its $\H$-Cholesky factorization. The advantage of the HDD method is that it does not require an agglomeration of the whole stiffness matrix. The memory is dynamically allocated and deallocated.

\begin{figure}[htbp]
\centerline{
\includegraphics[width=0.49\textwidth]{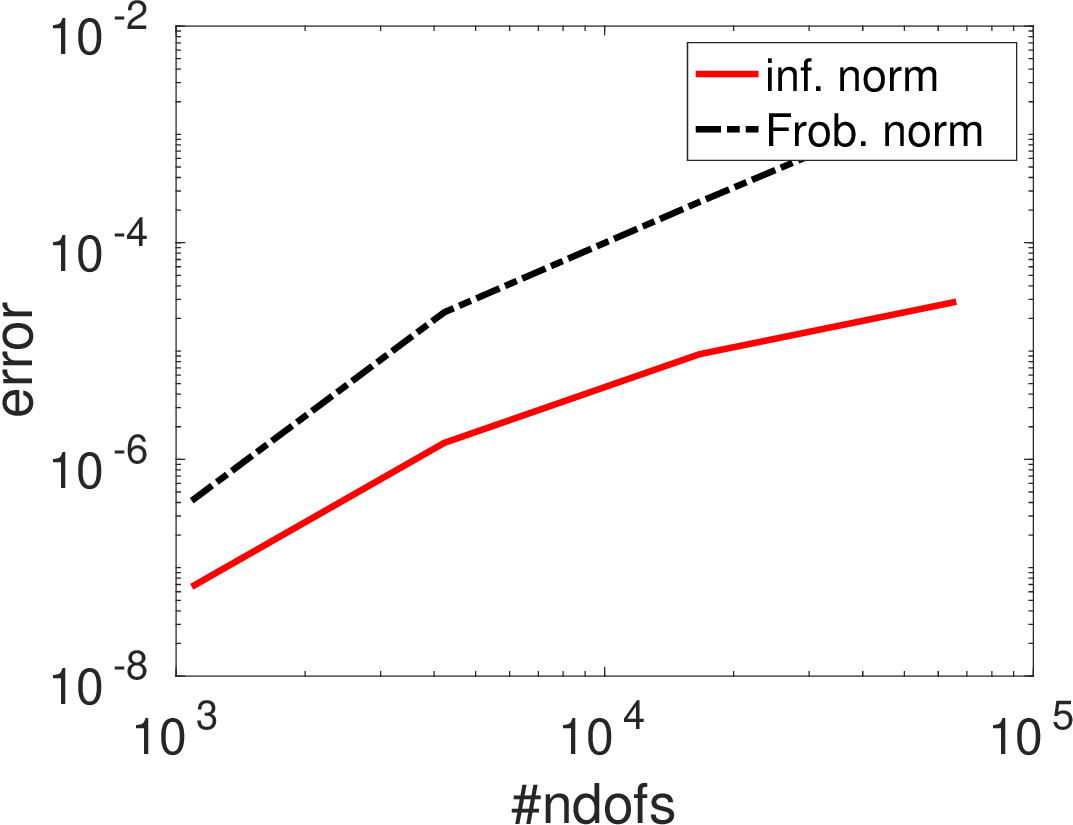}
\includegraphics[width=0.49\textwidth]{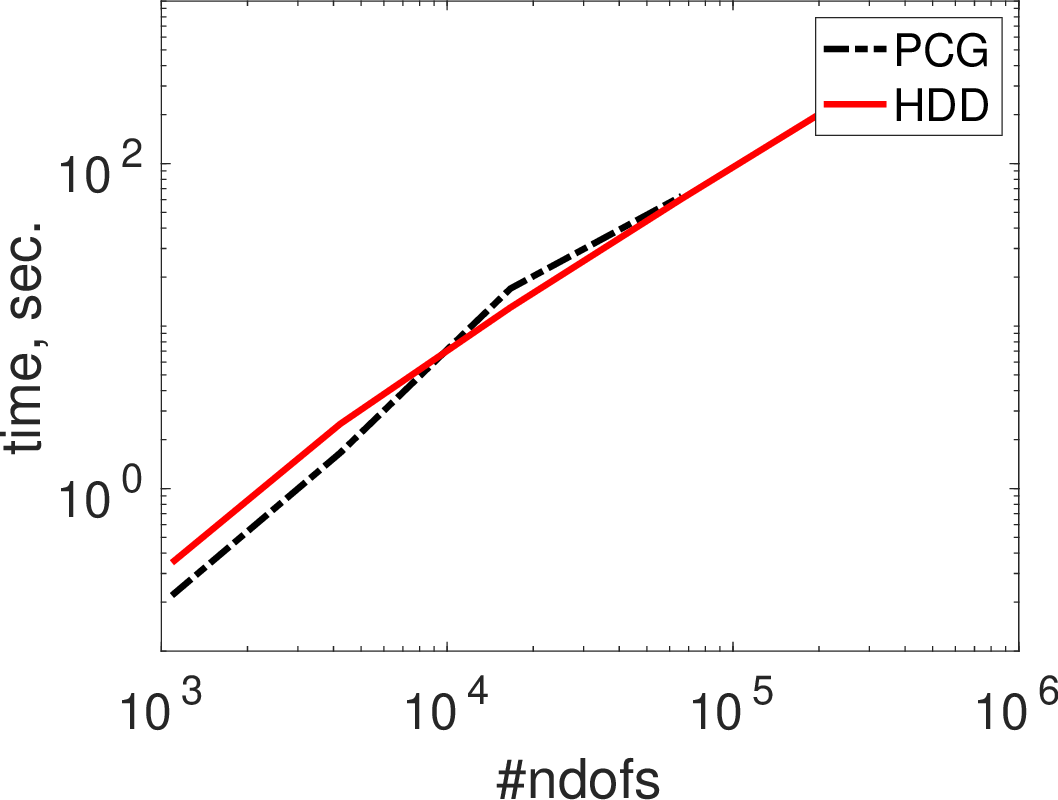}
}
\caption{Dependence of the absolute errors on the number of dofs, $f=1$, $\alpha(x,y)=1/(1.0001+\sin(500x)\sin(500y))$. $\H$-matrix accuracy $\varepsilon=10^{-5}$, $\frac{H}{h}=2$. (left) errors
(left) dependence of the errors $\Vert \tilde{\u}_{cg}-\tilde{\u}\Vert_2$ and $\Vert \tilde{\u}_{cg} - \tilde{\u} \Vert_{\infty}$ on $n$; (right) computing times PCG and HDD vs. $n$.}
\label{fig:hdd_pcg_sol_times}
\end{figure}

\begin{figure}[htbp]
\centerline{
\includegraphics[width=0.49\textwidth]{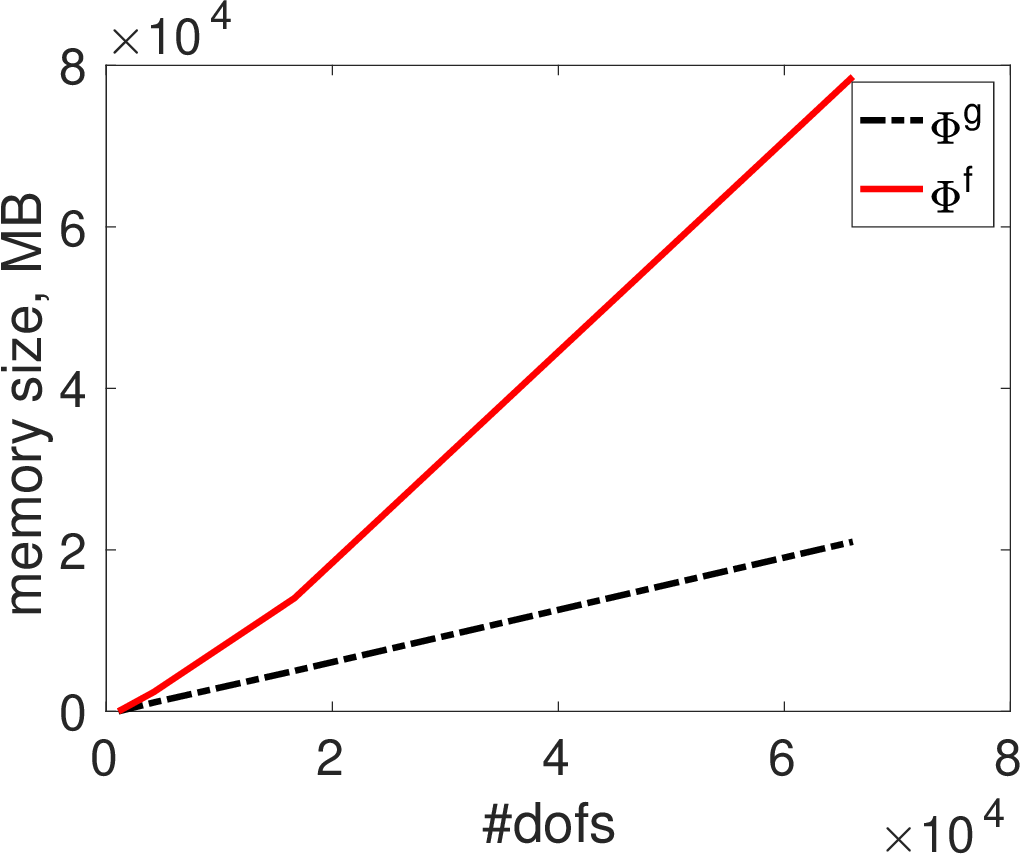}
}
\caption{Dependence of the total memory requirement for all $\Phi^g_{\omega}$ and $\Phi^f_{\omega}$ on $n$, the maximal $\H$-matrix rank is $k=7$.}
\label{fig:hdd_memory_phi}
\end{figure}

In the next example we take a coarse mesh for the right-hand side with the grid step size $H=2h$.
n Figure~\ref{fig:hdd_pcg_sol_times}(left) we visualise the difference $\Vert \tilde{\u}_{cg} - \tilde{\u} \Vert$ in the Frobenius and infinity norms. Here $\tilde{\u}$ is the HDD solution and $\tilde{\u}_{cg}$ the solution obtained by the PCG method with the $\H$-Cholesky preconditioner. 
The corresponding HDD and PCG times are compared in Fig.~\ref{fig:hdd_pcg_sol_times}(right). The accuracy inside of each $\H$-matrix subblock is $10^{-5}$. Note, that PCG requires too much memory for $n=513^2$ dofs and we were not able to compute $\tilde{\u}_{cg}$.\\ 
%
Figure~\ref{fig:hdd_memory_phi} shows the total storage requirement for all matrices $\Phi^g_{\omega}$ and $\Phi^f_{\omega}$, $\omega \in T_{\Th}$. We see an almost linear dependence on $n$.  

\section{Conclusion and discussion}
\label{sec:Conclusion}
We suggested another useful application of the already known HDD method. Namely, HDD speeds up computations of the data misfit (mismatch) in the likelihood function in the Bayesian approach.
HDD can also be used when the simulated data and measurement data are compared (e.g., in regression, parameter inference, data assimilation, Kalman filter, and Bayesian update problems).
HDD uses the fact that often only a functional of the solution or a small part of it is observed or measured. Therefore, HDD computes only a part of the inverse operator and only a part of the solution. Optimally, HDD computes only what is needed, i.e., what is measured.
 
As such the computational accuracy is as usual (for instance, as in the standard FEM method), but the computational recourses (FLOPS and storage) are smaller.
The HDD method is based on the hierarchical (recursive) domain decomposition, FEM, and the Schur complement methods. If the forward operator and the right-hand side can be discretized on different meshes, which is often the case in multiscale problems, the HDD method can get significant advantages. The computational resources will be reduced even more.

Additionally, to speed up the Schur complement computations, we approximate all intermediate and auxiliary matrices in the $\H$-matrix format. We then achieve the computational cost $\mathcal{O}(n\log^3 n)$ and the storage $\mathcal{O}(n\log^2 n)$. There is some overhead due to the construction of the hierarchical decomposition tree $T_{\Th}$ and permutation of indices.

To apply the HDD method, the user should have a possibility to 1) modify the assembling procedure of the stiffness matrix; 2) build the hierarchical domain decomposition tree. 

We note that the interface size in a $d$-dimensional problem is $\mathcal{O}(n^{d-1})$. So in a 2D case, the interface is $\mathcal{O}(n)$, whereas, in 3D problems, the interface is $\mathcal{O}(n^2)$. This fact results in increasing matrix sizes. The structure and the cost of the $\H$-matrix arithmetics become more expensive too.

Numerical tests showed that HDD requires more computational resources than PCG with the $\H$-Cholesky preconditioner, and less resources than the direct $\H$-matrix inverse. But the HDD method computes more details than PCG. It computes the solution operators $\Phi^g$ and $\Phi^f$ on each level of the hierarchical domain decomposition tree. These can be used later to compute a functional of the solution, or solution on different scales, in a sub-domain, in a point or on an interface.

The HDD method can be coupled with more uncertainty quantification and parameter inference techniques. Potentially interesting could be the coupling with the Multi-Level Monte Carlo method.

\subsection*{Acknowledgments}
This work was supported
by funding from the Alexander von Humboldt foundation (chair of Mathematics for Uncertainty Quantification at RWTH Aachen).

\newpage
\bibliographystyle{plain}
\bibliography{HDD_BU_inference} 

\newpage 
\begin{appendices}
\section{Appendix A}
\label{appendix:A}%
%

\begin{figure}[htbp]
\centerline{\includegraphics[width=0.99\textwidth]{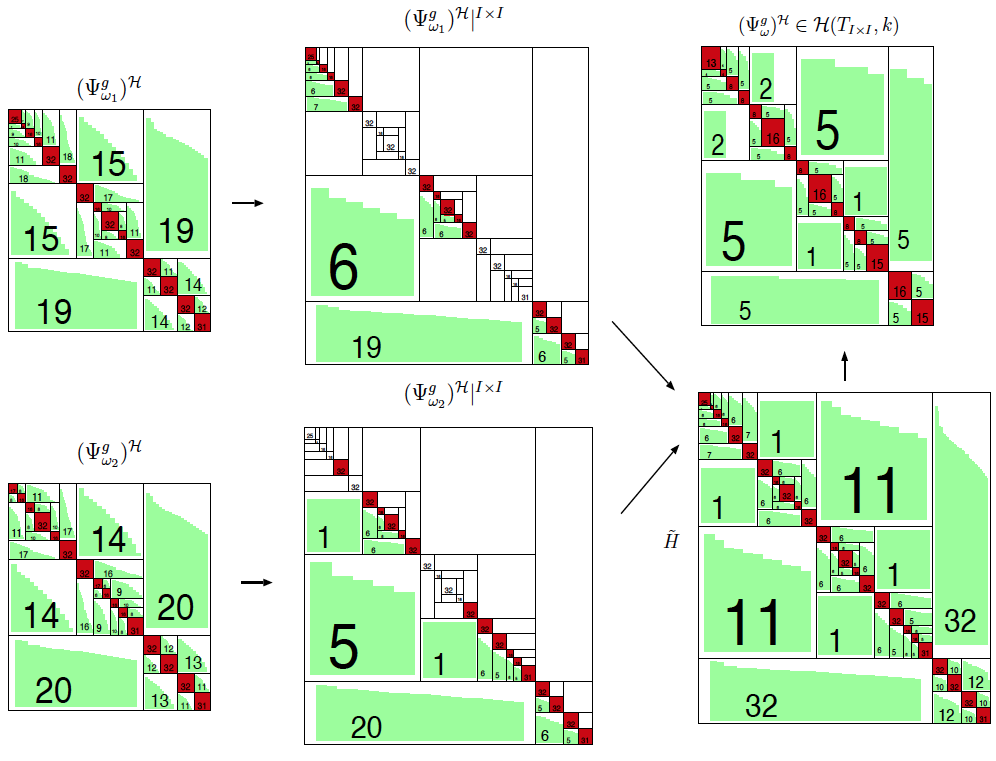}}
\caption{Building $(\Psi_{\omega}^g)^\H \in \mathbb{R}^{512 \times 512}$ from $(\Psi_{\omega_1}^g)^\H$ and $(\Psi_{\omega_2}^g)^\H \in \RR^{384 \times 384}$. The intermediate matrix $\Ht \in \mathbb{R}^{639 \times 639}$ is an auxiliary matrix. The maximal size of the diagonal blocks is $32 \times 32$. The red (dark) blocks indicate dense matrices. The green (grey) blocks indicate low-rank matrices. The steps inside these blocks show an exponential decay of the corresponding singular values. The white blocks indicate zero blocks. For the acceleration of building the symmetry of $\Psi_{\omega}^g$ is used.}
\label{fig:h2hg}
\end{figure}

\begin{ex}
Figure \ref{fig:h2hg} shows an example of building $(\Psi_{\omega}^g)^\H \in \mathbb{R}^{512 \times 512}$ from $(\Psi_{\omega_1}^g)^\H \in \RR^{384 \times 384}$ and $(\Psi_{\omega_2}^g)^\H\in \RR^{384 \times 384}$. 
Let $I:=I(\partial\omega\cup \gamma)$. The construction is performed in three steps:\\
1) embed matrix $(\Psi^g_{\omega_1})^{\H}$ into a larger matrix $H':=(\Psi^g_{\omega_1})^{\H}\vert^{I\times I}$ and $(\Psi^g_{\omega_2})^{\H}$ into $H'':=(\Psi^g_{\omega_2})^{\H}\vert^{I\times I}$,\\
2) since $H'$ and  $H''$ have the same $\H$-matrix format, compute the sum $\tilde{H}=H' \oplus H''$,\\
3) compute the Schur complement and eliminate the block (2,2) of size $I(\gamma) \times I(\gamma)$. \\
Note that $H'$, $H''$, $\Ht$ have the same block structures. The symmetries of $(\Psi^g_{\omega_1})^{\H}$, $(\Psi^g_{\omega_2})^{\H}$ and $(\Psi^g_{\omega})^{\H}$ are used. See more details and a similar construction of $(\Psi_{\omega}^f)^\H $ in \cite{MYPHD}.
\end{ex}
%
\end{appendices}
\end{document}